%% file: main.tex
\documentclass[11pt]{article}

\usepackage{xcolor}
\usepackage{geometry}   
\usepackage{setspace}
\usepackage{color}
\definecolor{clemson-orange}{RGB}{234,106,32}
\definecolor{chicago-maroon}{RGB}{128,0,0}
\definecolor{northwestern-purple}{RGB}{82,0,99}
\definecolor{cornell-red}{RGB}{179,27,27}
\definecolor{sauder-green}{RGB}{171,180,0}
\definecolor{harveymudd-gold}{RGB}{178,139,51}
\definecolor{lawngreen}{RGB}{0,250,154}
\definecolor{gray}{RGB}{192,192,192}
\definecolor{gray_header}{gray}{0.95}
\definecolor{purple_highlight}{HTML}{F3F0FF} 

\usepackage[small,compact]{titlesec}
\titlespacing*\section{0pt}{12pt plus 4pt minus 2pt}{4pt plus 0pt minus 0pt}
\titlespacing*\subsection{0pt}{6pt plus 4pt minus 2pt}{2pt plus 4pt minus 2pt}
\titlespacing*\subsubsection{0pt}{12pt plus 4pt minus 2pt}{0pt plus 2pt minus 2pt}

\usepackage[title]{appendix}
\usepackage{natbib}
 \bibpunct[, ]{(}{)}{,}{a}{}{,}%
 \setcitestyle{citesep={;}}

\usepackage{subcaption}
\usepackage{amssymb}
\usepackage{amsmath}
\usepackage{mathtools}
\usepackage{multicol}
\usepackage{graphicx}
\usepackage{array}
\usepackage{bm}
\usepackage{soul}
\usepackage{dirtytalk}
\usepackage{multirow}



\usepackage{algorithm}
\usepackage[noend]{algpseudocode}
\usepackage{caption}

\makeatletter
\newcommand\fs@boxedtopcap{\def\@fs@cfont{\bfseries}\let\@fs@capt\floatc@plain
	\def\@fs@pre{\setbox\@currbox\vbox{\hbadness10000
			\moveleft3.4pt\vbox{\advance\hsize by6.8pt
				\hrule \hbox to\hsize{\vrule\kern3pt
					\vbox{\kern3pt\box\@currbox\kern3pt}\kern3pt\vrule}\hrule}}}%
	\def\@fs@mid{\kern2pt}%
	\def\@fs@post{}\let\@fs@iftopcapt\iftrue}
\makeatother

\floatstyle{boxedtopcap}
\restylefloat{algorithm}

\usepackage{threeparttable}

\usepackage{enumitem}
\setlist[enumerate]{noitemsep, topsep=2pt}
\setlist[itemize]{noitemsep, topsep=2pt}

\usepackage{extarrows}

\usepackage{anysize}
\usepackage{latexsym}
\usepackage{amsthm}
\usepackage{epsfig}

\usepackage[T1]{fontenc}
\usepackage{mathrsfs}
\usepackage{float}
\usepackage{setspace}
\usepackage{color}
\definecolor{lawngreen}{RGB}{0,250,154}

\usepackage{titletoc}

\usepackage[colorlinks,citecolor=northwestern-purple,urlcolor=chicago-maroon,linkcolor=chicago-maroon]{hyperref}
\usepackage[nameinlink]{cleveref}
\crefname{assumption}{Assumption}{Assumptions}
\crefname{lemma}{Lemma}{Lemmas}
\crefname{theorem}{Theorem}{Theorems}
\crefname{corollary}{Corollary}{Corollaries}
\crefname{proposition}{Proposition}{Propositions}
\crefname{condition}{Condition}{Conditions}
\crefname{claim}{Claim}{Claims}
\crefname{procedure}{Procedure}{Procedures}
\crefname{algorithm}{Algorithm}{Algorithms}
\crefname{figure}{Figure}{Figures}
\crefname{remark}{Remark}{Remarks}
\crefname{section}{Section}{Sections}
\crefname{procedure}{Procedure}{Procedures}
\crefname{example}{Example}{Examples}
\crefname{definition}{Definition}{Definitions}
\crefname{table}{Table}{Tables}
\crefname{equation}{}{}
\crefname{enumi}{}{}
\crefname{conjecture}{Conjecture}{Conjectures}
\crefname{step}{Step}{Steps}
\crefname{appendix}{Appendix}{Appendices}
\crefname{footnote}{Footnote}{Footnotes}

\setlist[enumerate]{noitemsep, topsep=0pt}

\usepackage{thmtools}

\theoremstyle{definition}

\newtheorem{prob}{Problem}



\usepackage{xr}

\newcommand{\optfull}{OPT-Principled}
\newcommand{\optshort}{PriOPT}

\numberwithin{subcase}{case}

\numberwithin{subsubcase}{subcase}

\usepackage{natbib}
 \bibpunct[, ]{(}{)}{,}{a}{}{,}%
\providecommand{\keywords}[1]
{
  \small	
  \textbf{\textit{Keywords---}} #1
}

\definecolor{pink}{RGB}{255,0,128}

\usepackage{booktabs}

\usepackage{tabularx}

\usepackage{adjustbox}

\usepackage[table]{xcolor}

\usepackage{booktabs}
\usepackage{tabularx}
\usepackage{makecell} 

\usepackage{listings}
\usepackage[most]{tcolorbox}
\newtcolorbox{promptbox}[1]{
  enhanced,
  breakable,                
  colback=blue!3!white,       
  colframe=blue!35!black,     
  colbacktitle=blue!40!black, 
  title=\textbf{#1},       
  fonttitle=\bfseries\large,
  sharp corners,
  boxrule=0.8pt,
  left=5pt, right=5pt, top=5pt, bottom=5pt,
  parbox=false             
}

\newcommand{\code}[1]{\texttt{#1}}

\begin{document}

\title{From Soliloquy to Agora: Memory-Enhanced LLM Agents with Decentralized Debate for Optimization Modeling}

\author{
\textbf{
Jianghao Lin\textsuperscript{1,*},
Zi Ling\textsuperscript{2,*,\textdagger},
Chenyu Zhou\textsuperscript{1},
Tianyi Xu\textsuperscript{1},
Ruoqing Jiang\textsuperscript{4,\textdagger},
}\\
\textbf{
Zizhuo Wang\textsuperscript{3},
Dongdong Ge\textsuperscript{1}
}\\[0.5em]
\textsuperscript{1} Antai College of Economics and Management, Shanghai Jiao Tong University\\
\textsuperscript{2} University of Chicago Booth School of Business\\
\textsuperscript{3} The Chinese University of Hong Kong, Shenzhen (CUHK-Shenzhen)\\
\textsuperscript{4} School of Economics and Management, Tsinghua University\\[0.5em]
\textit{linjianghao@sjtu.edu.cn \quad zling@chicagobooth.edu \quad chenyuzhou@sjtu.edu.cn}\\
\textit{crimsonflag@sjtu.edu.cn \quad jiangrq@sem.tsinghua.edu.cn \quad wangzizhuo@cuhk.edu.cn}\\
\textit{ddge@sjtu.edu.cn}
}

\date{}

\begingroup
\renewcommand{\thefootnote}{\fnsymbol{footnote}}
\maketitle
\footnotetext[1]{Equal contribution.}
\footnotetext[2]{Corresponding authors.}
\endgroup
\setcounter{footnote}{0}

\begin{abstract}
\noindent Optimization modeling underpins real-world decision-making in logistics, manufacturing, energy, and public services, but reliably solving such problems from natural-language requirements remains challenging for current large language models (LLMs). 
In this paper, we propose \emph{Agora-Opt}, a modular agentic framework for optimization modeling that combines decentralized debate with a read-write memory bank. Agora-Opt allows multiple agent teams to independently produce end-to-end solutions and reconcile them through an outcome-grounded debate protocol, while memory stores solver-verified artifacts and past disagreement resolutions to support training-free improvement over time. This design is flexible across both backbones and methods: it reduces base-model lock-in, transfers across different LLM families, and can be layered onto existing pipelines with minimal coupling. 
Across public benchmarks, Agora-Opt achieves the strongest overall performance among all compared methods, outperforming strong zero-shot LLMs, training-centric approaches, and prior agentic baselines. Further analyses show robust gains across backbone choices and component variants, and demonstrate that decentralized debate offers a structural advantage over centralized selection by enabling agents to refine candidate solutions through interaction and even recover correct formulations when all initial candidates are flawed.
These results suggest that reliable optimization modeling benefits from combining collaborative cross-checking with reusable experience, and position Agora-Opt as a practical and extensible foundation for trustworthy optimization modeling assistance.
Our code and data are available at \url{https://github.com/CHIANGEL/Agora-Opt}.
\end{abstract}

\keywords{Large language models, optimization modeling, agentic debate, agentic memory}

\maketitle



\section{Introduction}

Operations research (OR) underpins decision-making in logistics, manufacturing, energy, and public services at a global scale \citep{singh2012overview,petropoulos2024operational}. At the center of these applications is \emph{optimization modeling}, which translates the operational challenges into mathematically well-posed decision variables, objectives, and constraints that deliver measurable impact in the real world. For instance, UPS’s ORION route-optimization system is reported to save about 10 million gallons of fuel annually and to generate \$300–\$400 million in yearly savings, alongside reductions of roughly 100{,}000 metric tons of CO$_2$ \citep{holland2017ups}. In another domain, the UN World Food Programme leveraged analytics and optimization to replan supply chains during COVID-19 and humanitarian crises, achieving more than \$150 million in savings while serving approximately 100 million people across over 80 countries \citep{peters2022world}. Despite such successes, building correct models directly from natural-language requirements remains a substantial obstacle for non-experts, and recent works on large language models (LLMs) have begun to narrow this gap by parsing problem text, producing formulations, and emitting solver-ready code \citep{ahmaditeshnizi2024optimus,huang2025orlm}.

Within this landscape, much of recent progress on applying LLMs to optimization modeling largely follows \emph{training-centric} approaches that update a base model, via fine-tuning based on instruction or reinforcement learning (RL), to improve the mapping from problem text to formulations and solver-ready code \citep{jiang2024llmopt,chen2025solver,huang2025orlm}. While effective, these approaches typically suffer from the \emph{base-LLM lock-in}: the trained model is tied to a specific base model version (e.g., Qwen2), so moving to a stronger successor does not transfer seamlessly. As a result, the substantial tuning invested in Qwen2 must often be repeated to obtain a trained model on Qwen2.5 or another base model.  

In parallel, \emph{agentic} methods have also been extensively studied because they treat the backbone as an interchangeable component and can adopt base-model upgrades with minimal additional adjustment. Classical agentic designs for optimization modeling instantiate a backbone LLM as one or more role agents to traverse the entire workflow.
For example, \citet{ahmaditeshnizi2024optimus} coordinate a team of agents (i.e., manager, formulator, programmer, and evaluator) to process structured optimization tasks. However, because backbones in these typical agentic methods are closed-source, these methods \emph{cannot directly benefit from training data} within the method. Once deployed, agents tend to operate as fixed systems that do not incorporate online experience. One commonly used technique, Retrieval-Augmented Generation (RAG), can ground responses in externally retrieved documents, but it remains read-only: even when seeded with existing data, the system cannot accumulate its own trial-and-error or solver-verified fixes unless the corpus is explicitly rebuilt, because there is no experience write-back.

A further limitation cuts across both families: most existing work, whether training-centric or agentic, performs reasoning at a given stage with a single backbone, exhibiting \emph{single-model myopia} that limits diversity and internal cross-checking. Even when an agentic framework invokes different backbones across steps, each step is usually executed by one backbone at a time, preserving model-specific idiosyncrasies but reducing robustness. One natural remedy is to introduce \emph{agentic debate}, which integrates the intelligence of multiple backbones to improve overall performance.  Debate frameworks generally fall into two types. In \emph{centralized} debate, a small set of backbones exchange arguments under a moderator that evaluates and aggregates the outcome \citep{liang2024encouraging,long2024multi}. While helpful, this setup inherits the biases and idiosyncrasies of the judge and thus does not fully resolve the robustness and diversity concerns above. In \emph{decentralized} debate, there is no moderator and decisions follow external tests or agreement among independently produced results \citep{chen2025debatecoder,li2025swe}. However, in practice, it can be challenging to obtain convergence, and different backbones may fail to agree or to achieve a better response without a well-defined consensus principle, raising questions about when to stop, how to reconcile near-ties, and how to arbitrate conflicting outputs.

\begin{figure*}[t]
  \centering
  \includegraphics[width=1\linewidth]{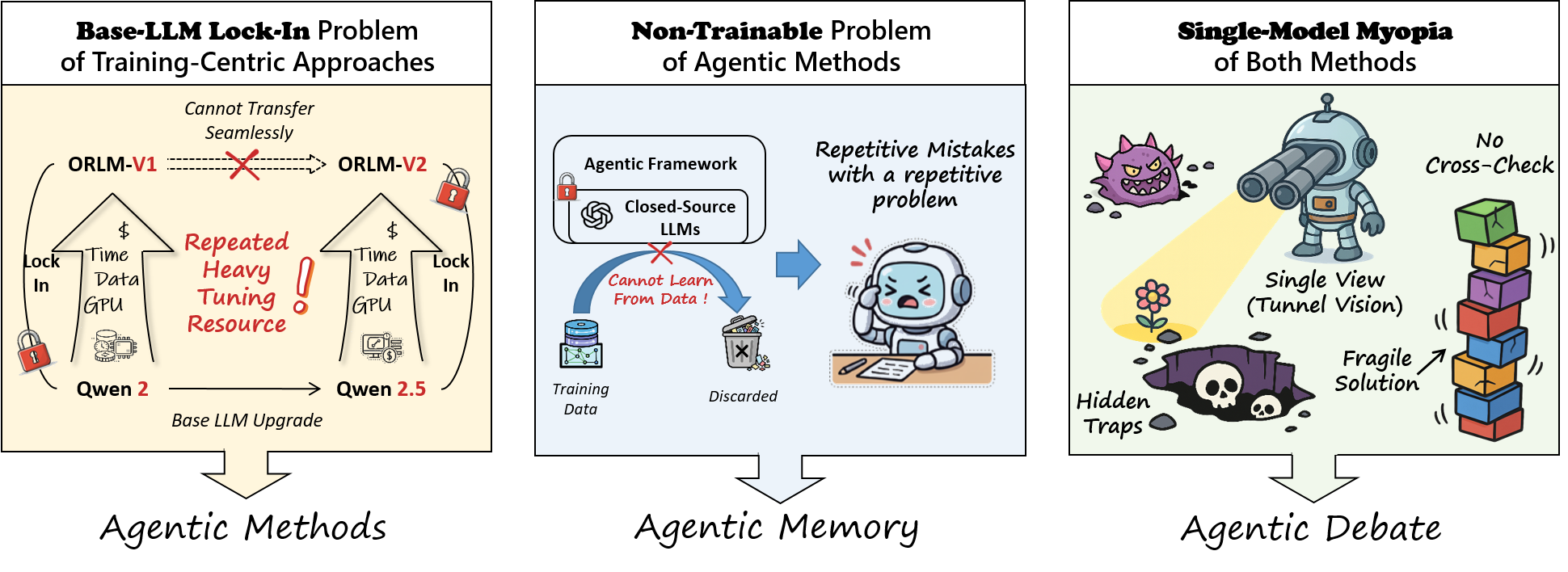}
  \caption{The illustration of three limitations in most existing methods: (a) base–LLM lock–in of training–centric approaches, (b) non–trainable problem of agentic methods, and (c) single–model myopia; alongside their paired design principles in our framework for LLM–based optimization modeling: an agentic foundation for easy backbone upgrades, a read–write agentic memory design, and decentralized agentic debate.}
  \label{fig:limitations}
\end{figure*}

As a brief guide to what follows, \cref{fig:limitations} synthesizes the three limitations discussed above and previews the targeted design principles of our framework that we will introduce next.

To address these limitations, we present \textbf{Agora-Opt}, a \emph{unified} agentic framework that couples decentralized \underline{\textbf{ag}}entic debate with a read–write agentic mem\underline{\textbf{or}}y b\underline{\textbf{a}}nk for \underline{\textbf{opt}}imization modeling. We begin from an \emph{agentic} foundation: the backbone LLM is treated as an interchangeable component within a role-structured pipeline that moves from problem text to formulation, code, and solver output. Consequently, when migrating to a new backbone, we reuse the same roles, procedures, and writeback loop without any internal parameter retuning or redesign. The same process applies across models, yielding strong mobility across backbones (evaluated in \cref{subsec:swap_backbone}).

Building on this agentic foundation, Agora-Opt tackles \emph{single-model myopia} through \emph{agentic debate}. We deliberately avoid a \emph{centralized} scheme with a moderator, because such judges inherit the biases and idiosyncrasies of their own backbone, and ultimately re-concentrate authority in a single model at adjudication time. Instead, we leverage a core characteristic of optimization modeling: although the solving pipeline is essentially multi-stage, the process culminates in \emph{quantitative endpoints} (solution and objective values) that can be checked independently of any judge. Accordingly, Agora-Opt adopts \emph{decentralized agentic debate}: multiple agent teams, instantiated on diverse backbones and methods, run in end-to-end manner, and the system outputs an answer only when their \emph{final, solver-verified} outcomes align or when the maximum number of debate rounds is reached. The consensus is thus defined by objective signals rather than subjective summaries, thereby better integrating diverse intelligence. In effect, this design avoids the moderator bias, highlights cross-backbone and cross-method discrepancies on challenging cases, and makes adjudication a measurable, solver-grounded criterion.

Agora-Opt then addresses the fixed-at-deploy behavior observed in prior agentic systems through a \emph{read–write agentic memory bank}. Beyond simple storage, the memory is designed to work with the debate mechanism. 
It comprises two complementary stores: a \textit{generation memory} that accumulates verified problem-solving experiences (e.g., formulating, coding, and debugging), and a \textit{debate memory} that preserves the argumentative reasoning and consensus patterns derived from multi-agent collaboration. This pairing is intentional: generation memory accelerates routine episodes to solve each problem, while debate memory preserves how disagreements were resolved, including what teams proposed, which checks mattered, and which fixes led to convergence, so future debates start from a richer base of verified experience. This memory bank also improves upgrade robustness, since it preserves solver-verified know-how across backbone changes and reduces the need for retraining or prompt retuning. Finally, because the optimization community and industry are facing new tasks and problems at a rapid pace, the memory enables online leverage of the real-time solving process, allowing the agent to improve between runs rather than only between model releases.

Finally, the framework is intentionally \emph{modular} and \emph{flexible across both backbones and methods}: decentralized debate operates on solver-verified endpoints and the memory writes back experience rather than tuning parameters, \emph{Agora-Opt} can be applied to a wide range of backbones and readily layered onto existing pipelines with minimal coupling. 
We specify our concrete team configurations in \cref{sec:methodology}, and in \cref{sec:analysis_discussion} we demonstrate flexibility of our proposed framework by instantiating teams with different backbones (e.g., Gemini, GPT, DeepSeek) and embedding representative methods (e.g., ORLM in \cite{huang2025orlm}, StepORLM in \cite{zhou2025steporlm}, OptiMUS in \cite{ahmaditeshnizi2024optimus}) to verify both effectiveness and portability.

In summary, our main contributions can be summarized as follows: 
\begin{enumerate}[label=(\roman*)]
\item \textbf{Framework and novelty.} We introduce \emph{Agora-Opt}, a unified agentic framework that, to the best of our knowledge, is the first to couple decentralized debate with an agentic memory for optimization modeling.
The design is modular and flexible across both backbones and methods, enabling \emph{Agora-Opt} to layer onto most existing pipelines and improve their performance while simultaneously mitigating base-LLM lock-in and reducing re-tuning cost when upgrading to stronger base models.

\item \textbf{Decentralized debate protocol.}  To our knowledge, we formalize the first debate protocol tailored to \emph{optimization modeling}. In this outcome-grounded scheme, multiple agent teams independently produce end-to-end solutions, and the system outputs only on \emph{consensus}. This removes single-model myopia, enables cross-checking, and combines collective intelligence across models and methods into a principled, quantitative convergence rule.

\item \textbf{Agentic memory design.} We develop a memory bank with \emph{generation} and \emph{debate} memories that write back per-task artifacts and outcomes, as well as how disagreements are resolved, tightly integrated with the debate process. This yields training-free improvement after deployment, preserves solver-verified know-how across backbone upgrades, and leverages rapidly released new tasks by integrating new experience without parameter updates.

\item \textbf{Evaluation and generality.} To evaluate our frameworks, we conduct extensive experiments on six public benchmarks together with \optfull, a curated testbed of challenging optimization instances derived from a public resource, comparing \emph{Agora-Opt} with strong zero-shot LLMs, training-centric and agentic baselines in \cref{sec:main_results}, showing that Agora-Opt achieves the most competitive overall performance among all compared methods. Beyond the main results, \cref{sec:analysis_discussion} studies the robustness and generality of the framework through backbone-swap experiments, component ablations, sensitivity analyses on debate rounds, and further demonstrates breadth by layering our protocol over existing methods such as ORLM, StepORLM, and OptiMUS. Crucially, we isolate the structural advantage of decentralized debate over centralized selection, demonstrating that interactive debate can repair and synthesize correct formulations even when all initial candidate solutions are flawed.
\end{enumerate}

The rest of the paper is organized as follows: \cref{sec:literature review} reviews related literature. \cref{sec:methodology} formally presents the Agora-Opt framework, including the agent-team design, the decentralized debate protocol, and the read-write memory mechanism. \cref{sec:main_results} reports the main experimental results on public OR benchmarks, comparing Agora-Opt against zero-shot LLMs, training-centric models, and agentic baselines, and includes a representative case study showing how decentralized debate and memory retrieval resolve modeling ambiguity and guide convergence to a final solution. \cref{sec:analysis_discussion} provides further analyses of Agora-Opt’s robustness and generality, including compatibility across different backbone LLMs, ablation studies of key components, and an in-depth examination of the decentralized debate protocol through comparisons with centralized judge selection, sensitivity analyses on debate rounds, and generalization behavior.
Finally, \cref{sec:conclusion} concludes the paper.

\section{Literature Review}
\label{sec:literature review}

Our paper is mainly related to three research streams: LLMs for operations research problems, agentic debate, and memory augmentation.

\vspace{2mm}
\noindent\textbf{LLMs for operations research problems.}
Large language models (LLMs) have recently been extensively explored to bridge the gap between natural-language problem descriptions and mathematical optimization models, producing solver-ready code. The NL4Opt competition offered an early demonstration and a widely used benchmark, showing that general-purpose LLMs can extract entities and structure from text to produce mathematical modeling formulations \citep{ramamonjison2023nl4opt}. Building on this foundation, subsequent systems couple LLMs with classical solvers to automate more of the pipeline from problem text to both formulation and code.

Following this trajectory, much recent progress adopts \emph{training-centric} approaches that update a base LLM with curated synthetic data and instruction/RL fine-tuning to internalize optimization-modeling knowledge.
On the fine-tuning side, LLMOPT \citep{jiang2024llmopt} employs multi-instruction tuning with alignment and self-correction to improve formulation and code generation, while Solver-Informed RL \citep{chen2025solver} grounds learning in verifiable solver feedback to reduce hallucinations and improve factual correctness. Recent advances in LLM training, including comparison-oracle preference learning and stepwise guided policy optimization to correct intermediate reasoning, further strengthen fine-tuning pipelines \citep{chen2025compo, chen2025stepwise}.
On the data synthesis side, ORLM \citep{huang2025orlm} introduces a two-stage data synthesis pipeline linking natural-language problems, formal formulations, and executable code before fine-tuning open-soured LLMs for end-to-end modeling. Building on this line of work, \cite{zhou2025auto} tailor DualReflect-style generation for dynamic programming. OptMath \citep{lu2025optmath} provides scalable bidirectional synthesis with forward modeling and rejection sampling, while Step-Opt \citep{wu2025step} increases task difficulty through iterative synthesis with structured validation. Complementing these, StepORLM \citep{zhou2025steporlm} further studies a self-evolving, process-supervised training scheme for OR language models. Despite these gains, training-centric systems inherit \emph{base-LLM lock-in}: trained models are tied to specific base releases, therefore upgrading to stronger base models often requires substantial re-tuning and fails to transfer seamlessly. 

An alternative line of work employs \emph{agentic methods}, treating a backbone LLM as one or more role agents and viewing the backbone as an interchangeable component that can adopt base-LLM upgrades with minimal adjustment.
Multi-agent designs illustrate this trend: the Chain-of-Experts framework proposed by \cite{xiao2023chain} coordinates terminology, modeling, programming, and reflection under a conductor agent. Building on this, \cite{ahmaditeshnizi2024optimus} introduce OptiMUS, which also uses a conductor agent to coordinate multiple steps, while refining each conversational step before dispatching tasks to the next agent. To mitigate conductor-driven unpredictability, \cite{wang2025ormind} propose a structured, cognitive-inspired workflow with counterfactual reasoning, named ORMind, to enhance the reliability and clarity of solutions.
Even with these advances, essential gaps remain: at any given stage, reasoning is typically executed by a single backbone (\emph{single-model myopia}). Moreover, since the backbone is closed-source, the method cannot directly benefit from training data within the agentic workflow, and long-horizon brittleness persists as recurrent specification and coding errors reappear across tasks.

\vspace{2mm}
\noindent\textbf{Agentic debate} has emerged as a robust paradigm to enhance the reasoning, factuality, and reliability of LLMs by moving beyond single-model systems, which are often constrained by their internal knowledge and fixed inferential patterns, to a collaborative system where multiple models interact to solve a problem.

The idea traces back to AI safety: \citet{irving2018ai} propose a self-play workflow in which two AI agents act as debaters to persuade a human judge. Their “arguments” are not natural-language statements but selections of individual image pixels to ``convince'' a basic classifier in a simple adversarial game, thereby illustrating the debate concept's potential.
With the advent of powerful LLMs, a series of works shifted from AI safety to performance improvements in reasoning and factual accuracy, retaining this “third-party” judge to control outcomes. For example, \citet{liang2024encouraging} introduce ``Multi-Agent Debate'' (MAD) framework to solve a key failure mode of single-LLM self-reflection: ``Degeneration-of-Thought'' (DoT), which means once confident in an incorrect solution, the model is unable to generate novel or divergent thoughts to correct itself. Subsequently, more researchers further investigate on how to optimize agentic debate's workflow \citep{long2024multi, estornell2024acc}, which makes the process more robust and efficient. Meanwhile, the purpose of the debate mechanism becomes more diverse, being applied not only for inference-time reasoning but also for auxiliary tasks such as evaluation \citep{chan2023chateval} and model training \citep{subramaniam2024debategpt}.
However, in all these papers, a third party is always required as a moderator, critic, or referee team to evaluate and summarize competing arguments, which we term as a \emph{centralized} debate. Although effective, centralized setups inherit third-party bias and errors, and can perpetuate \emph{single-model myopia} when the judge’s backbone dominates decisions. 

To address these issues, some recent work turns to \emph{decentralized} debate, where outcomes are decided by objective checks or by agreement among independently produced results rather than a judge.
\citet{du2023improving} designed a framework employing multiple identical LLM instances that propose, exchange, and iteratively converge through debate without a separate judging model. \citet{liu2025breaking} introduce Diverse Multi-Agent Debate (DMAD), allowing agents to follow diverse reasoning paths and collectively arrive at an answer. Recent applications of agentic debate have moved into structured, high-stakes domains, such as software engineering (SWE-Debate \citep{li2025swe}) and code generation (DebateCoder \citep{chen2025debatecoder}), but have not yet entered operations research.
Distinct from these prior lines, we introduce \emph{decentralized} debate into optimization modeling, leveraging its inherently quantitative endpoints to adjudicate agreement across diverse backbones.

\vspace{2mm}
\noindent\textbf{Memory augmentation} equips LLMs with external stores that expand or persist knowledge beyond fixed parameters.
The canonical method is Retrieval-Augmented Generation (RAG), firstly proposed by \citet{lewis2020retrieval}, which grounds generation in external corpora via a retriever–generator architecture \citep{wu2024retrieval}. While improving factuality, RAG is read-only mainly and static over largely external/static information (docs, manuals), and cannot accumulate the agent’s own experiences unless the corpus is thoroughly re-built.

Recent research has shifted toward dynamic \emph{read-write} memory architectures that support continual agent evolution \citep{modarressi2023ret}. These mechanisms have been widely applied in social simulation \citep{park2023generative}, system-level memory management \citep{packer2023memgpt, xu2025mem}, and long-term interaction handling via forgetting mechanisms \citep{zhong2024memorybank}. In highly structured domains, memory architectures evolve to store verified skills and logical insights rather than unstructured text: \citet{wang2023voyager} utilize a ``skill library'' for executable code, while \citet{wang2025mirix} maintain memory for long-horizon reasoning. 
This trajectory is further formalized by \citet{zhou2025memento}, who demonstrate that such memory-based updates can effectively substitute for parameter fine-tuning to enable continuous agent evolution.
Specifically in operations research, \citet{kong2025alphaopt} construct an evolving experience library to store structured insights, covering both domain modeling and solver syntax, by refining their applicability conditions over time.

Distinct from these works, we design a more comprehensive and flexible memory architecture for \emph{Agora-Opt} that comprises both \emph{generation memory} and \emph{debate memory}. While \citet{kong2025alphaopt} manage a hierarchical taxonomy of abstracted insights tailored for operations research, such a structured design can limit generalization: rigid taxonomic boundaries often struggle to capture the structural dependencies and knowledge transfer across diverse problem domains.
In contrast, our framework supports flexible write-in and read-out access to \textit{generation memory} and, uniquely, incorporates a \emph{debate memory} that preserves entire episodes of the consensus-building process.
By recording the complete trajectory of how diverse agents reconcile disagreements, we enable the system to reuse not only verified solutions but also effective collaboration strategies for resolving ambiguity in multi-agent scenarios.

\input{texts/method}

\input{texts/4-results}

\input{texts/5-analysis}

\section{Conclusion}
\label{sec:conclusion}

In this paper, we propose \emph{Agora-Opt}, a memory-enhanced agentic framework with decentralized debate for optimization modeling. The framework treats backbone LLMs as interchangeable, leverages cross-backbone diversity through debate, and improves through read--write experience accumulation without retraining. To the best of our knowledge, Agora-Opt is the first framework to combine decentralized agentic debate with a read--write memory bank for optimization modeling.

Empirically, across six public benchmarks together with \textsc{OPT-Principled}, Agora-Opt achieves the strongest overall performance among all compared methods. Our analyses further show that its gains are robust across backbone choices and component variations, and that decentralized debate offers a clear advantage over centralized judge selection by enabling refinement through interaction and recovering correct answers even when both initial teams are wrong. In addition, Agora-Opt can be layered over representative existing pipelines, highlighting its modularity and generality.

These gains reinforce our central point: independent teams surface discrepant assumptions, debate turns these discrepancies into targeted, execution-grounded refinements, and memory converts resolved failures into reusable know-how, echoing how humans iteratively collaborate across diverse perspectives while scaling beyond individual expertise.
Looking forward, two promising directions are: (i) \emph{adaptive control of debate}, which dynamically determines when debate is needed and how much is required based on disagreement and observed progress; and (ii) \emph{extension to harder and more realistic OR tasks}, where our findings suggest particularly strong benefits in high-difficulty regimes. Together, these directions position Agora-Opt as a practical and extensible foundation for trustworthy optimization modeling assistance that can evolve with rapidly changing tasks and backbone models.

\bibliographystyle{plainnat}
\bibliography{custom}

\newpage
\section*{\centering \LARGE Online Appendix}
\vspace{1.5em}  


\begin{appendices}

\input{appendix/benchmarks.tex}
\input{appendix/prompts.tex}

\end{appendices}

\end{document}

%% file: texts/method.tex
\section{Methodology}
\label{sec:methodology}



\begin{figure*}[t]
  \centering
  \includegraphics[width=1\linewidth]{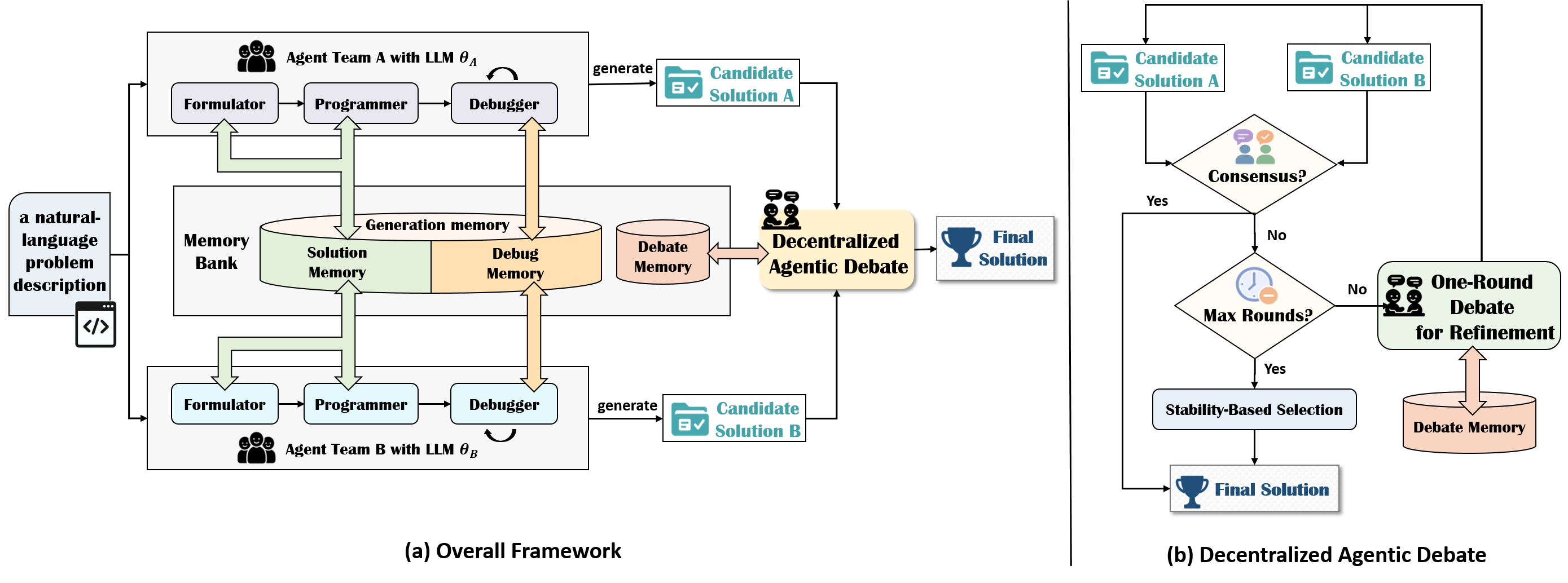}
  \caption{Overview of the Agora-Opt framework.
    \textbf{(a) Overall framework.}
    A natural language optimization problem is solved by two symmetric agent teams (Formulator–Programmer–Debugger) built on different backbone LLMs, which interact with a unified memory bank and feed their candidate solutions into a decentralized agentic debate.
    \textbf{(b) Decentralized agentic debate.}
    The two candidate solutions enter a debate: if they reach consensus, a final solution is returned; otherwise, guided by debate memory, the system either performs one more refinement round or applies stability-based selection when the round budget is exhausted.}
  \label{fig:overview}
\end{figure*}

In this section, we formally present our agentic framework Agora-Opt for optimization modeling.
As illustrated in \cref{fig:overview}, Agora-Opt takes a natural-language description of an optimization problem as input and routes it to two agent teams that share the same internal roles, prompts and workflows, but are built on different backbone LLMs. Each team follows a three-stage pipeline: it first formulates the problem, then generates solver code, and finally executes and debugs the code to obtain a candidate solution, along with its objective value and diagnostic logs. The two candidate solutions are then processed by an agentic debate protocol, which consists of a trigger mechanism, an iterative refinement loop and a termination condition. Throughout both single-team solving and multi-team debate, a unified memory bank provides read and write augmented experience, where solution memory and debug memory support formulating, programming, and debugging inside each team, and debate memory stores past reconciliation traces that guide future debates. Accordingly, our framework is organized into three key components: \emph{the agent-team generation}, \emph{agentic debate protocol}, and \emph{agentic memory design}. The full set of prompts used for all agent roles and stages is provided in Appendix~\ref{appendix_prompts}.

\subsection{The Agent-team Generation}
\label{subsec:agent_team_generation}

Agora-Opt adopts a fixed agent-team design paired with different backbone LLMs. 
Unless otherwise stated, our main framework deploys two such agent teams that share the same roles, prompts, and workflow, differing only in the underlying LLM.
In \cref{subsec:swap_methods}, we further replace this default agent team with alternative methods, including other agentic designs and training-centric models, to assess the generalization of our framework.

Formally, we denote the natural-language problem description as $x \in \mathcal{X}$. An agent team, parametrized by a backbone LLM $\theta$, defines a mapping $\mathcal{T}_\theta: \mathcal{X} \to \mathcal{S}$, where $\mathcal{S}$ denotes the space of candidate solutions. For a given problem $x$, the agent team generates a candidate solution $s \in \mathcal{S}$ as a tuple $s = (f, c, v, \mathcal{L})$, where:
\begin{itemize}
    \item $f$ is the structured mathematical formulation;
    \item $c$ is the executable solver code;
    \item $v \in \mathbb{R} \cup \{\bot\}$ represents the solver-evaluated objective value, where $v \in \mathbb{R}$ indicates a successful execution and $\bot$ denotes a failure;
    \item $\mathcal{L}$ denotes the \textit{execution feedback}, containing solver logs, error tracebacks, and warning messages required for diagnosis.
\end{itemize}
The generation process proceeds in three sequential stages: \emph{formulating}, \emph{programming}, and \emph{executing and debugging}.
Conceptually, these stages are realized as three specialized modules (illustrated as Formulator, Programmer, and Debugger), all implemented via LLM calls.

\paragraph{Formulating.} The formulator $\Phi_{\text{form}}$ first parses the optimization task and then maps the input $x$ to a structured formulation:
\begin{equation}
    f = \Phi_{\text{form}}(x; \theta).
\end{equation}
This formulation specifies the decision variables, objective function, and constraints. Optionally, this process can be augmented by retrieving solution patterns from the memory bank (as detailed in \cref{subsec:agentic_memory}) to guide the formulating.

\paragraph{Programming.} Conditioned on the formulation $f$, the programmer $\Phi_{\text{prog}}$ translates the structured specification into executable solver code (e.g., Gurobi-based Python):
\begin{equation}
    c = \Phi_{\text{prog}}(x, f; \theta)
\end{equation}
Similar to the formulating stage, code fragments retrieved from the memory bank can optionally serve as in-context references to ensure syntax correctness, which is introduced in \cref{subsec:agentic_memory}.

\paragraph{Executing and Debugging.} The debugger $\Phi_{\text{debug}}$ executes the generated code $c$ within a solver environment $\mathcal{E}$, which returns an outcome tuple:
\begin{equation}
    (v, \mathcal{L}) \leftarrow \mathcal{E}(c).
\end{equation}
If the solver terminates successfully ($v \neq \bot$), the agent team accepts the run and extracts a candidate solution $s=(f,c,v,\mathcal{L})$.
However, if execution fails ($v = \bot$), then the agent team enters a debugging cycle: the debugger inspects error messages and solver logs $\mathcal{L}$ to edit the code $c$ for a corrected revision, optionally querying the debug memory for historical fix strategies (detailed in \cref{subsec:agentic_memory}):
\begin{equation}
    c \leftarrow \Phi_{\text{debug}}(x, f, c; \theta).
\end{equation}
The revised code $c$ is then re-executed.
This execute–debug–revise loop iterates until either a successful execution is obtained ($v \neq \bot$) or a predefined retry budget $K_{\text{retry}}$ is exhausted.

At the end of this process, each agent team outputs a candidate solution $s$ for instance $x$, including its formulation $f$, code $c$, and execution outcomes $(v, \mathcal{L})$; the latter $(v, \mathcal{L})$ may indicate success or contain error diagnostics if the retry budget was exhausted, allowing the subsequent debate process to further revise the formulation and code. 
These candidates are then passed to the agentic debate protocol.

\subsection{Agentic Debate Protocol}
\label{subsec:agentic_debate}

After the above agent-team generation stage, Agora-Opt typically has two candidate solutions per instance, produced by two agent teams with the same workflow design but different backbone LLMs. 
Let $\mathcal{T}_A$ and $\mathcal{T}_B$ denote two agent teams parametrized by distinct backbones $\theta_A$ and $\theta_B$. Given an input $x$, they independently produce initial solutions $s_A^{(0)}$ and $s_B^{(0)}$. 
The agentic debate protocol consists of three components: a \textit{trigger mechanism} that initiates the debate, an \textit{iterative refinement} loop for solution improvement, and a \textit{consensus and termination} criterion to finalize the output.

\paragraph{Trigger Mechanism.} The debate is not activated indiscriminately. It is triggered only when the initial solutions exhibit a substantive disagreement, defined as either a \textit{feasibility discrepancy} (either one fails) or an \textit{optimality gap} (both succeed but yield different values). Formally, the trigger condition is:
\begin{equation}
    \left(v_A^{(0)} = \bot\right) \; \lor \; \left(v_B^{(0)} = \bot\right) \; \lor \; \left(\left|v_A^{(0)} - v_B^{(0)}\right| > \epsilon\right),
    \label{eq:consensus condition}
\end{equation}
where $\epsilon$ is a predefined tolerance threshold. If this condition holds, then the system initiates the debate process; otherwise, it proceeds directly to the phase of \textit{consensus and termination}.

\paragraph{Iterative Refinement.} 
In each debate round $t\ge 1$, both teams engage in a peer-review process: each team examines the problem description $x$, together with the current pair of its own and the opponent's solutions $\left(s_A^{(t-1)}, s_B^{(t-1)}\right)$. 
They are prompted to identify formulation mistakes, missing or redundant constraints, and inconsistencies of the objective functions in both solutions, and then propose revised formulations and codes.
The revised codes are re-executed, and the new outcomes $\left(s_A^{(t)}, s_B^{(t)}\right)$ are used in the next round.
Crucially, both teams have the option to query the debate memory (see \cref{subsec:agentic_memory}) for summarized precedents on how similar disagreement patterns have been resolved.
At the end of each round $t$, a branching condition governs the workflow: if the maximum round limit is reached ($t = T_{\text{max}}$), the protocol transitions directly to the termination phase; otherwise, the updated solutions are re-evaluated against the trigger criteria in \cref{eq:consensus condition} to determine whether a subsequent round of debate is needed or not.

\paragraph{Consensus and Termination.} 
The debate terminates when the trigger condition \cref{eq:consensus condition} no longer holds or the maximum round budget $T_{\text{max}}$ is exhausted. 
If both candidates execute successfully and their objective values converge within a tolerance $\epsilon$, we treat this as a consensus and thereby achieve a final \emph{debate-refined} solution, typically the better-performing of the two.
If convergence is not reached by $t=T_{\text{max}}$, we apply a stability-based fallback: we select the candidate that changes the least between the last two rounds, under the intuition that a more stable team is more confident in its formulation after debate.

\subsection{Agentic Memory Design}
\label{subsec:agentic_memory}

All stages of Agora-Opt are supported by a unified memory bank $\mathcal{M}$ storing reusable experience from past tasks. 
As outlined in the introduction, this bank consists conceptually of a \emph{generation memory} and a \emph{debate memory}. 
Motivated by our agent-team design, we further decompose the generation memory into two components:
\emph{solution memory}, which stores successful problem–formulation–code triplets, and \emph{debug memory}, which stores failure–repair episodes.
Hence, $\mathcal{M}$ is composed of three disjoint sets: solution memory $\mathcal{M}_{\text{sol}}$, debug memory $\mathcal{M}_{\text{bug}}$, and debate memory $\mathcal{M}_{\text{deb}}$. 
We employ a dense vector-based retrieval mechanism based on a semantic embedding model $E(\cdot)$.

\paragraph{Retrieval Function.} We define a generic retrieval operator $\mathcal{R}$ for querying any memory component $\mathcal{M}_* = \{(k_n,\text{val}_n)\}_{n\in[|\mathcal{M}_*|]}$, where $k_n$ represents the semantic key used for indexing, and $\text{val}_n$ denotes the stored artifact. Given a query vector $q = E(\text{query\_content})$, the operator returns the top-$N$ most relevant entries based on cosine similarity:
\begin{equation}
    \mathcal{R}(q,\mathcal{M}_*, N)
= \operatorname*{arg\,max}_{\substack{S \subseteq \mathcal{M}_* \\ |S| = N}}
    \sum_{(k_n,\text{val}_n) \in S} \cos\bigl(q, E(k_n)\bigr).
\end{equation}
Specifically, we implement $E(\cdot)$ based on a public text embedding model \texttt{bge-small-en-v1.5}\footnote{\url{https://huggingface.co/BAAI/bge-small-en-v1.5}}.

\paragraph{Solution Memory.} Solution memory stores successful problem-formulation-code triplets that have passed execution checks:
\begin{equation}
    \mathcal{M}_{\text{sol}} = \{ (x_n, (f_n, c_n)) \mid v_n \neq \bot \}.
\end{equation}
Here, the problem description $x_n$ serves as the retrieval key. During the formulating and programming stages for a new problem $x'$, the agent retrieves $\mathcal{R}(E(x'), \mathcal{M}_{\text{sol}},n)$ to inject proven formulation patterns and code syntax into the prompt context to guide the new solution generation.

\paragraph{Debug Memory.} 

This component is constructed from cases in which initially generated code failed but was subsequently repaired by the executing-and-debugging loop.
Each entry is indexed by an error signature $sig(\mathcal{L})$, which is an LLM-generated minimal problem context (summarizing the question, formulation, and code). Formally, debug memory is defined as:
\begin{equation}
    \mathcal{M}_{\text{bug}} = \{ (sig(\mathcal{L}_n), (\mathcal{L}_n, diagnosis_n, fix_n)) \},
\end{equation}
where $diagnosis_n$ is the LLM’s explanation of the root cause, and $fix_n$ is the corresponding fix strategy or corrected code sketch.  
When a new execution fails with a solver log $\mathcal{L}'$, the debugger queries $\mathcal{R}(E(sig(\mathcal{L}')), \mathcal{M}_{\text{bug}},N)$ to retrieve the diagnosis and fix strategy summarized from similar past failures, thereby guiding the code revision.

\paragraph{Debate Memory.} 
Debate memory is built from debate runs in which two agent teams started with a substantial disagreement and later converged to a consensus solution:
\begin{equation}
    \mathcal{M}_{\text{deb}} = \{ (\text{concat}(x_n, \Delta_n), \mathcal{H}_{\text{deb},n}) \}.
\end{equation}
The retrieval key is the concatenation of the problem $x_n$ and the initial discrepancy description $\Delta_n$ (i.e., LLM-summarized conflicts between candidate solutions). 
The stored value $\mathcal{H}_{\text{deb},n}$ contains the key arguments exchanged during debate (such as pointing out missing constraints or mis-specified objectives),  and the final consensus formulation, optionally accompanied by an LLM-written summary of the decisive evidence and recommended formulating pattern.
The debate memory allows agents to reuse reconciliation experiences and logical checks for resolving similar ambiguities, leading to more effective \emph{debate strategies} for optimization modeling.

%% file: texts/4-results.tex
\section{Main Results}
\label{sec:main_results}

\subsection{Experiment Setups}


\subsubsection{Benchmarks.}
To evaluate Agora-Opt on optimization modeling and solving across varying difficulty levels and problem types, we conduct our evaluation on six diverse public benchmarks that are widely used in the OR community: NL4Opt~\citep{Ramamonjison2023}, MAMO (split into EasyLP and ComplexLP)~\citep{huang2024mamo}, NLP4LP~\citep{AhmadiTeshnizi2024}, ComplexOR~\citep{xiao2024chainofexperts}, IndustryOR~\citep{huang2025orlm} and ReSocratic~\citep{yang2025optibenchmeetsresocraticmeasure}. More details on these six public benchmarks are provided in Appendix~\ref{appendix:benchmarks}.

We further include \textsc{OPT-Principled}, a curated benchmark of challenging optimization instances derived from the public \textsc{OPT-Engine} framework~\citep{chen2026optenginebenchmarkinglimitsllms}. \textsc{OPT-Engine} programmatically generates solver-verifiable optimization problems with controllable mathematical and semantic complexity, and its analysis shows that LLM difficulty is strongly shaped by both mathematical scale and non-canonical constraint structure. Guided by these principles, we construct \textsc{OPT-Principled} by selecting instances that emphasize these challenging characteristics and better reflect the scale and constraint richness of realistic OR tasks.

\subsubsection{Baselines.}
To ensure a comprehensive comparison, we compare against a dirverse set of representative methods spanning three categories:

\begin{itemize}
    \item \textbf{Zero-shot LLMs:} We evaluate leading general-purpose LLMs in a zero-shot setting, including \texttt{OpenAI-o3}~\citep{openai_o3_2025}, \texttt{Gemini-2.5-Pro}~\citep{comanici2025gemini}, \texttt{GPT-4o}~\citep{openai2024gpt4o}, \texttt{Kimi-K2}~\citep{team2025kimi}, \texttt{DeepSeek-R1}~\citep{guo2025deepseek}, \texttt{DeepSeek-V3}~\citep{deepseekai2025deepseekv3technicalreport}, \texttt{Qwen2.5-72B-Instruct}~\citep{qwen2025qwen25technicalreport}, \texttt{Qwen3-32B}~\citep{yang2025qwen3technicalreport}, and \texttt{Qwen3-8B}~\citep{yang2025qwen3technicalreport}.
    \item \textbf{Training-centric Models:} We include state-of-the-art fine-tuned models specifically designed for OR tasks, including \texttt{ORLM}~\citep{huang2025orlm}, \texttt{LLMOPT}~\citep{jiang2024llmopt}, \texttt{OptMATH}~\citep{lu2025optmath}, SIRL~\citep{chen2025solverinformedrlgroundinglarge}, and \texttt{StepORLM}~\citep{zhou2025steporlm}.
    \item \textbf{Agentic Methods:} We compare against recent agentic frameworks tailored for OR tasks, including \texttt{OptiMUS}~\citep{AhmadiTeshnizi2024}, Chain-of-Experts (\texttt{CoE})~\citep{xiao2024chainofexperts}, Chain-of-Thought (\texttt{CoT}) \citep{wei2022chainofthought}, and \texttt{CAFA}~\citep{deng24cafa}.
\end{itemize}

\subsubsection{Implementation Details and Evaluation Metrics.}

We implement Agora-Opt in Python and, unless otherwise specified, our primary experiments instantiate the two heterogeneous agent teams with \texttt{GPT-4o} and \texttt{DeepSeek-V3} as the backbone LLMs.
For all generation calls, we set the maximum context length to 16,384 tokens to support long optimization problem descriptions and memory contents, and we use a temperature of $T=0.01$ by default to ensure reproducibility and stability. 
Beyond single-team solving, our standard setting enables a debate protocol that is triggered when the two teams’ solver-verified outcomes disagree (predefined tolerance threshold $\epsilon=5\times 10^{-2}$) and is capped at 3 rounds (\(T_\text{max}=3\)).
Under this configuration, the unified memory bank is queried via vector similarity: we retrieve the top-$N$ ($N=4$) entries from \textit{solution memory}, the top-$N$ ($N=3$) entries from \textit{debug memory} upon execution failures, and the top-$N$ ($N=2$) entries from \textit{debate memory} to guide reconciliation. We use Gurobi as the solver for all generalist and agentic method evaluations, and we additionally enable the executing-and-debugging loop with up to 3 retries per instance and a 120-second execution timeout during the generation phase.
For fair comparison, we implement the agentic baselines using both \texttt{GPT-4o} and \texttt{Gemini-2.5-Pro} as backbone LLMs. 
Since the open-source code of \texttt{CoE} is tightly coupled with a specific version of LangChain and \texttt{GPT-4o}, we are unable to seamlessly switch the backbone to \texttt{Gemini-2.5-Pro}; as a result, we only report \texttt{CoE}'s performance using \texttt{GPT-4o}.

Finally, for evaluation, we mark an instance as correct if the generated solution's objective value matches the ground truth within a relative tolerance of $5\%$ ($\varepsilon=0.05$), switching to an absolute tolerance of $10^{-3}$ when the ground-truth objective is zero (since zero is not dividable). The evaluation execution timeout is set to 90 seconds.

\subsection{Overall Performance}

\input{tables/main_results}

We report the comparative results of Agora-Opt against state-of-the-art baselines across six public operations research benchmarks together with \textsc{OPT-Principled} in \cref{table:main_results}. 
Overall, our proposed Agora-Opt achieves a new state-of-the-art macro-average Pass@1 accuracy of 84.6\%, significantly outperforming both the strongest training-centric model (\texttt{StepORLM}, 71.4\%) and the most capable zero-shot LLM (\texttt{OpenAI-o3}, 75.6\%).

Beyond the headline average, this table reveals a consistent narrative: Agora-Opt achieves \emph{experience-augmented collective intelligence} by coordinating heterogeneous models through debate while continually reusing experience from memory, leading to more reliable optimization modeling across benchmarks.
We highlight four key observations from \cref{table:main_results} regarding the efficacy of agentic debate, the lock-in dilemma of training-centric methods, the data utilization gap in prior agentic methods, and the advantage of our framework on harder benchmarks.

\textbf{Emergent Intelligence from Weaker Backbones.}
A central finding is that Agora-Opt achieves superior performance by synthesizing backbones that are individually weaker than those of the top-tier proprietary models.
As shown in the ``Zero-shot LLMs'' category, the two base LLMs used by Agora-Opt, \texttt{GPT-4o} (63.9\%) and \texttt{DeepSeek-V3} (67.9\%), lag behind frontier models such as \texttt{OpenAI-o3} (75.6\%) or \texttt{Gemini-2.5-Pro} (74.3\%).
However, by integrating them via the Agora-Opt framework, the system achieves a substantial performance boost (+16.7\% over \texttt{DeepSeek-V3} and +20.7\% over \texttt{GPT-4o} on average), ultimately surpassing even \texttt{OpenAI-o3}.
This validates that our decentralized debate mechanism and read-write memory can effectively mitigate single-model myopia, enabling the whole system to exceed the capabilities of its constituent parts without relying on the most expensive frontier models.

\textbf{Addressing LLM Lock-in in Training-Centric Methods.}
The results further underscore the rapid depreciation of training-centric methods due to base-model lock-in.
While fine-tuning offers strong performance at release, these models struggle to keep pace with the rapid evolution of general-purpose LLMs.
For instance, \texttt{ORLM} and \texttt{LLMOPT}, which were competitive upon release, now average around 58.6\%, falling behind modern zero-shot models like \texttt{DeepSeek-V3} (67.9\%) and \texttt{Kimi-K2} (69.6\%).
Even the current SOTA training-centric model, \texttt{StepORLM} (71.4\%), maintains only a slim margin over standard inference models like \texttt{OpenAI-o3} (75.6\%) and \texttt{Gemini-2.5-Pro} (74.3\%).
This suggests a diminishing return on fine-tuning investments: as base models improve, the ``sunk cost'' of training previous versions becomes a liability.
In contrast, Agora-Opt is agnostic to the backbone and can instantly leverage stronger models without retraining.

\textbf{Bridging the Data Utilization Gap in Agentic Systems.}
Prior agentic methods (e.g., \texttt{OptiMUS}, \texttt{CoE}) avoid training lock-in, but they often face a pronounced ``data utilization gap''.
While these methods can strengthen LLM-based optimization solving by introducing structured and multi-step workflows, they typically still fall short of data-rich training-centric models like \texttt{StepORLM}, largely because they operate as ``read-only'' systems that cannot internalize verified optimization data.
Agora-Opt bridges this gap by utilizing a read-write agentic memory, which preserves the flexibility of agentic inference while continually retrieving verified formulation patterns, debugging traces, and debate trajectories. As a result, the agentic framework Agora-Opt can compound useful experience over time and break through the performance ceiling of heavy training-centric approaches.

\textbf{Robust Performance in the Hard Regime.}
The advantage of Agora-Opt is more pronounced in the challenging benchmarks, which are historically brittle for LLM-based optimization modeling and where small formulation mistakes, missing constraints, or execution failures occur more frequently.
For instance, on \texttt{IndustryOR}, where strong baselines remain in the 50–70\% range (e.g., \texttt{GPT-4o}: 51.0\%, \texttt{DeepSeek-V3}: 56.0\%, \texttt{OpenAI-o3}: 70.0\%), Agora-Opt reaches 74.0\%, yielding large gains over its constituent backbones (+23.0\% and +18.0\% over \texttt{GPT-4o} and \texttt{DeepSeek-V3} respectively). This “hard-tail lift” is further reflected by \texttt{MAMO}’s difficulty split: while single-model baselines often collapse from \texttt{EasyLP} to \texttt{ComplexLP} (e.g., \texttt{DeepSeek-V3}: 95.2\%$\rightarrow$53.2\%, \texttt{OpenAI-o3}: 93.9\%$\rightarrow$63.1\%), Agora-Opt shows a much smaller degradation (97.6\%$\rightarrow$81.1\%). This advantage is also evident on \textsc{OPT-Principled}, where strong zero-shot baselines reach only 49.0\%--57.0\%, the strongest prior agentic baseline reaches 64.0\%, and Agora-Opt attains 70.0\%. Taken together, these results indicate that Agora-Opt scales more gracefully with instance complexity and semantic non-canonicity, making it a promising foundation for future OR tasks where problem statements are longer, constraints are more intricate, and solver/code interactions become increasingly error-sensitive.

\subsection{Case Study}
\label{subsec:case_study}

To explicitly demonstrate the mechanics of Agora-Opt, and in particular how agentic memory and decentralized debate jointly shape the evolution of the two agent teams' solutions, we present a qualitative analysis of a representative instance: the \textit{Paint Mixing Problem} (Problem 398 in benchmark ReSocratic).

\begin{figure*}[t]
  \centering
  \includegraphics[width=1\linewidth]{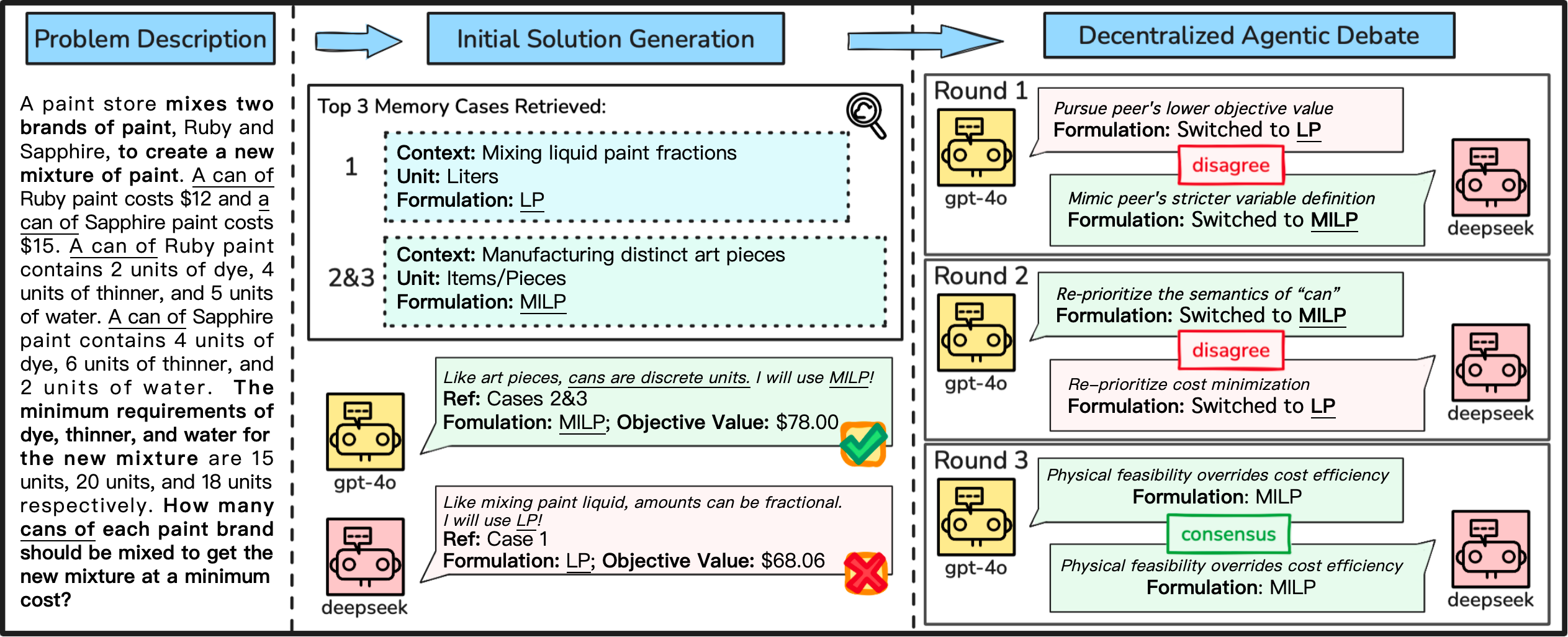}
  \caption{An illustrative trace of Agora-Opt solving the Paint Mixing Problem via agentic memory and decentralized debate}
  \label{fig:case_study}
\end{figure*}

As illustrated in \cref{fig:case_study}, this task is to minimize the cost of mixing two paint brands while satisfying chemical component requirements.
The core modeling challenge arises from semantic ambiguity: although the chemical components (dye, thinner) naturally suggest continuous decision variables, the purchasing unit of ``cans'' implicitly requires discreteness. However, the problem statement does not explicitly specify integrality, leaving the choice between Linear Programming (LP) and Mixed-Integer Linear Programming (MILP) open to interpretation.

\paragraph{Contradictory Retrieval References}
In the \textit{Initial Solution Generation} phase, both agents retrieve the same top-$3$ historical cases from the solution memory. As described in \cref{subsec:agentic_memory}, the system attempts to inject proven formulation patterns from these retrieved examples. However, in this specific instance, the retrieval introduces \textbf{contradictory modeling patterns}:
\begin{itemize}
    \item \textbf{The Semantic Distractor (Case 1):} A paint blending problem featuring a \emph{continuous pattern} (LP), where variables represent liquid fractions.
    \item \textbf{The Structural Analogs (Cases 2 \& 3):} Art production problems featuring a \emph{discrete pattern} (MILP), where variables represent indivisible units.
\end{itemize}
Given these mixed signals, \texttt{DeepSeek-V3} (Agent Red) places more weight on the surface-level semantic match to “paint” in Case 1 and adopts a relaxed LP formulation. In contrast, GPT-4o (Agent Yellow) prioritizes the deeper structural match to indivisible units in Cases 2 \& 3, leading it to the intended MILP formulation.


\paragraph{Dynamics of Decentralized Agentic Debate.}
The subsequent debate process highlights the dynamics of our decentralized protocol, which relies not on a central moderator or fixed stopping rules, but on the spontaneous convergence of independent agents.
In each round, both agents are prompted to revise their previous generation by cross-referencing their own formulation with the peer's proposal and execution outcomes.
This iterative synthesis drives the following trajectory:

\begin{enumerate}
    \item \textbf{Round 1 (The Crossover):} A phenomenon we term the \emph{Objective Value Trap} emerges. Agent Yellow, observing Agent Red's significantly lower cost (\$68.06 vs. \$78.00), relaxes its integrality constraints to pursue the apparent efficiency gain. Simultaneously, Agent Red, critically evaluating Agent Yellow's stricter variable definition, adopts the MILP formulation. This leads to a complete role reversal, where both agents swap formulations but remain in disagreement.
    \item \textbf{Round 2 (Oscillation and Correction):} The debate acts as a stress test for constraint hardness. Agent Yellow re-examines the physical semantics of ``cans'' and returns to the MILP formulation, prioritizing feasibility over optimality. In contrast, Agent Red reverts to the LP formulation, seemingly pulled back by the mathematical allure of the lower objective value.
    \item \textbf{Round 3 (Consensus):} Continued critique leads to convergence. Notably, Agent Yellow maintains its MILP formulation from the previous round, while Agent Red ultimately accepts the integrality requirement implied by the purchasing unit and aligns with it. Both teams therefore output the integer-feasible solution (\$78.00), demonstrating that physical feasibility acts as a hard constraint that must override nominal cost improvements from.
\end{enumerate}
This trajectory clarifies the complementary roles of memory and debate in Agora-Opt. Solution memory expands the hypothesis space by providing diverse (even conflicting) patterns, whereas decentralized debate acts as a crucial verification filter, effectively pruning the attractive but invalid relaxed formulation that a single model might otherwise accept.

%% file: tables/main_results.tex
\begin{table*}[t]
\footnotesize
\centering
\caption{The overall performance of Agora-Opt and baselines with Pass@1 accuracy (\%) on six OR benchmarks.
Best results are highlighted in \textbf{bold} and the second-highest values are \underline{underlined}.
Since \texttt{OptMATH} is not publicly available, we only report the scores cited from its original publication, marked with the symbol (*), while missing entries are denoted with (-).}
\label{table:main_results}

\setlength{\tabcolsep}{2pt}
\adjustbox{max width=1.0\textwidth}{
\begin{tabularx}{\textwidth}{l*{10}{>{\centering\arraybackslash}X}}

\toprule
{\scriptsize Model} & {\scriptsize Params} & {\scriptsize NL4OPT} & \multicolumn{2}{c}{\scriptsize MAMO} & {\scriptsize NLP4LP} & {\scriptsize CompOR} & {\scriptsize IndOR} & {\scriptsize ReSoc} & {\scriptsize \optshort} & {\scriptsize Avg.} \\
\cmidrule(lr){4-5}
 & & & {\scriptsize EasyLP} & {\scriptsize ComplexLP} & & & & & & \\
\midrule

\rowcolor{gray!15}
\multicolumn{10}{c}{\textbf{\textit{Zero-shot LLMs}}} \\

\texttt{OpenAI-o3} & Closed
& 78.4
& 93.9
& \textbf{63.1}
& \underline{93.8}
& \textbf{72.2}
& \textbf{70.0}
& 84.4
& \underline{49.0}
& \textbf{75.6}
\\

\texttt{Gemini-2.5-Pro} & Closed
& \textbf{82.6}
& 87.9
& 52.3
& \textbf{94.9}
& \underline{66.7}
& 63.0
& \textbf{89.6}
& \textbf{57.0}
& \underline{74.3}
\\

\texttt{GPT-4o} & Closed
& 76.1
& 93.6
& 54.1
& 91.6
& 44.4
& 51.0
& 84.1
& 16.0
& 63.9
\\

\texttt{Kimi-K2} & 1T
& 77.9
& 93.4
& 55.9
& 84.3
& 61.1
& 55.0
& 81.9
& 47.0
& 69.6
\\

\texttt{DeepSeek-R1} & 671B
& 77.5
& 90.3
& \underline{59.5}
& 87.6
& \underline{66.7}
& \underline{67.0}
& 83.9
& 47.0
& 72.4
\\

\texttt{DeepSeek-V3} & 671B
& \underline{79.8}
& \underline{95.2}
& 53.2
& 92.1
& 55.6
& 56.0
& \underline{85.1}
& 26.0
& 67.9
\\

\texttt{Qwen2.5-72B-Instruct} & 72B
& 78.9
& \textbf{95.8}
& 44.1
& 88.2
& 50.0
& 41.0
& 81.1
& 29.0
& 63.5
\\

\texttt{Qwen3-32B} & 32B
& 77.5
& 92.3
& 46.9
& \underline{93.8}
& 50.0
& 54.0
& \underline{85.1}
& 47.0
& 68.3
\\

\texttt{Qwen3-8B} & 8B
& 63.8
& 73.6
& 45.0
& 58.4
& 27.8
& 48.0
& 61.0
& 32.0
& 51.2
\\

\midrule
\midrule

\rowcolor{gray!15}
\multicolumn{10}{c}{\textbf{\textit{Training-centric Methods}}} \\

\texttt{ORLM} & 8B
& 73.8
& 90.4
& 59.5
& 76.4
& \underline{50.0}
& 38.0
& 61.8
& 19.0
& 58.6 \\
  
\texttt{LLMOPT} & 14B
& 75.1
& 83.5
& 67.6
& 86.0
& 22.2
& \underline{39.0}
& 73.2
& \textbf{22.0}
& 58.6
\\

\texttt{OptMATH} (origin) & 32B
& \underline{95.9*}
& 89.9*
& 54.1*
& - & - & - & - & - & - \\

\texttt{SIRL} & 7B
& 94.8
& \textbf{97.3}
& \textbf{82.0}
& \underline{93.3}
& \underline{50.0}
& 38.0
& \underline{80.1}
& \underline{21.0}
& \underline{69.6}
\\

\texttt{StepORLM} & 8B
& \textbf{97.7}
& \underline{97.2}
& \underline{79.3}
& \textbf{97.8}
& \textbf{55.6}
& \textbf{40.0}
& \textbf{82.6}
& \underline{21.0}
& \textbf{71.4}
\\

\midrule
\midrule

\rowcolor{gray!15}
\multicolumn{10}{c}{\textbf{\textit{Agentic Methods}}} \\

\texttt{OptiMUS-v0.3(GPT-4o)} & Closed
& 76.2
& 78.0
& 46.8 
& 88.8
& 44.1
& 30.0
& 87.6
& 20.0
& 58.9
\\

\texttt{OptiMUS-v0.3(Gemini)} & Closed
& \underline{80.3}
& 82.0
& 54.1
& 92.7
& 61.1
& 44.0
& 88.8
& \underline{64.0}
& 70.9
\\

\texttt{CoT(GPT-4o)} 
& Closed 
& 62.2 
& 49.5 
& 42.3 
& 74.7 
& 39.2 
& 36.0
& 43.6 
& 31.0
& 47.3
\\

\texttt{CoT(Gemini)} 
& Closed 
& 72.4
& 91.4
& 45.0
& 89.3
& 50.0
& \underline{61.0}
& 86.6
& 35.0
& 66.3 
\\

\texttt{CAFA(GPT-4o)} 
& Closed 
& 68.1 
& 71.2 
& 44.5 
& 50.0 
& 46.4 
& 39.0
& 40.1
& 24.0
& 47.9
\\

\texttt{CAFA(Gemini)} 
& Closed 
& \textbf{92.5}
& \underline{94.5}
& \textbf{82.0} 
& \underline{97.8}
& \underline{61.1}
& 52.0
& \underline{96.0}
& 62.0
& \underline{79.7}
\\ 

\texttt{CoE(GPT-4o)} 
& Closed 
& 66.7 
& \underline{94.5} 
& 50.6
& 87.4 
& 57.1 
& 45.0
& 71.2 
& 18.0
& 61.3
\\

\rowcolor{purple!15}
\texttt{Agora-Opt (GPT-4o+DS-V3)} & Closed & \textbf{92.5} & \textbf{97.6} & \underline{81.1} & \textbf{98.9} & \textbf{66.7} & \textbf{74.0} & \textbf{96.3} & \textbf{70.0} & \textbf{84.6} \\

\bottomrule
\end{tabularx}
}
\par
\scriptsize
\textbf{Abbreviations}: CompOR: ComplexOR, IndOR: IndustryOR,  \optshort: \optfull, Avg: Macro-Average, DS-V3: DeepSeek-V3.
\end{table*}

%% file: texts/5-analysis.tex
\section{Analysis and Discussion}
\label{sec:analysis_discussion}

\subsection{Compatibility w.r.t. Backbone LLMs}
\label{subsec:swap_backbone}

To verify that the effectiveness of \texttt{Agora-Opt} is not tied to a specific base model (i.e., mitigating the ``lock-in'' issue discussed in \cref{table:main_results}), we evaluate our framework across different backbone combinations.
As shown in \cref{table:agora_variants}, we compare the corresponding single-agent baselines (\texttt{Gemini-2.5-Pro}, \texttt{DeepSeek-V3}, \texttt{GPT-4o}) with three Agora-Opt variants instantiated from these models.
Overall, the results reveal three critical insights regarding model interoperability, the value of heterogeneous debate, and robustness on high-difficulty instances.

\input{tables/LLM_compatibility}

\textbf{Model-agnostic Enhancement.}
First of all, \texttt{Agora-Opt} consistently outperforms its constituent single models regardless of the backbone configuration.
All three paired variants achieve a macro-average accuracy between 84.6\% and 85.4\%, significantly surpassing the best single-model baseline (\texttt{Gemini-2.5-Pro} at 74.2\%).
Notably, even when pairing the relatively weaker \texttt{GPT-4o} with \texttt{DeepSeek-V3}, the system still reaches 84.6\%, effectively narrowing the gap to proprietary frontier models.
This consistency confirms that our framework's performance gains arise from its unified process of reasoning, debate, and memory-based experience reuse, rather than the raw capability of any single underlying LLM.

\textbf{The Potential of Diversity Beyond Raw Strength.}
Beyond being model-agnostic, the results suggest that the diversity of the agent team may be more critical than the individual strength of the agents.
Surprisingly, although \texttt{GPT-4o} has the lowest single-agent score (63.9\%), the overall performance of \texttt{GPT-4o + Gemini} (85.3\%) remains nearly identical to that of \texttt{DS-V3 + Gemini} (85.4\%), suggesting that a weaker individual backbone does not necessarily lead to a weaker collaborative system. Moreover, the combination of \texttt{GPT-4o + Gemini} achieves stronger results on several benchmarks, including \texttt{MAMO-ComplexLP} (84.7\% vs.\ 82.0\%), \texttt{NLP4LP} (99.4\% vs.\ 98.9\%), and \texttt{IndustryOR} (83.0\% vs.\ 76.0\%).
Its inclusion appears to contribute different reasoning patterns or error distributions that complement \texttt{Gemini}'s capabilities more effectively than \texttt{DeepSeek-V3} on these tasks.
This provides supporting evidence for the ``collective intelligence'' hypothesis: a diverse team of heterogeneous agents can cross-verify and correct each other more efficiently than a homogeneous team of stronger but potentially correlated experts.

\textbf{Robustness in High-Difficulty Regimes.}The synergistic effect is most pronounced on the hardest benchmarks, aligning with our earlier observation in \cref{table:main_results}.
On \texttt{IndustryOR}, where the strongest single model (\texttt{Gemini}) scores only 63.0\% and \texttt{GPT-4o} struggles at 51.0\%, their combination in \texttt{Agora-Opt} leaps to 83.0\%, representing a remarkable 20\% absolute improvement over the best individual component.
A similar pattern also appears on \texttt{\optshort}: although the improvement is less dramatic, paired systems still substantially outperform their constituent single models, with \texttt{DS-V3 + Gemini} reaching 75.0\% compared with 57.0\% for \texttt{Gemini} and 26.0\% for \texttt{DeepSeek-V3}. 
Taken together, this indicates that \texttt{Agora-Opt} is especially valuable in high-difficulty regimes, where the debate mechanism, supported by memory retrieval, effectively leverages the complementary strengths of different backbones to solve complex, ambiguity-laden problems that any single model struggles to handle in isolation.

\subsection{Ablation Study}

To further verify the contribution of each component in our framework, we conducted an ablation study on three representative challenging benchmarks, \texttt{ComplexLP}, \texttt{IndustryOR}, and our developed \texttt{\optfull}, where modeling ambiguity and failure modes are more prevalent.
As shown in \cref{table:ablation_study}, we start with the full \texttt{Agora-Opt} framework and sequentially remove the debate memory, debug memory, solution memory, and the debate mechanism, ultimately reducing the framework to a vanilla LLM.
The results demonstrate that every module contributes to the final performance, though their impact varies across benchmarks.

\input{tables/ablation_study}

The results show that while the agentic structure provides a strong baseline, the memory components are the primary drivers of state-of-the-art performance.
Removing the \textbf{Solution Memory} causes the sharpest decline on \texttt{ComplexLP} (from 76.6\% to 64.9\%), indicating that retrieving formulation templates is crucial for logically complex problems.
Conversely, the \textbf{Debug Memory}, proves essential for \texttt{IndustryOR} (from 72.0\% to 64.0\%) and \texttt{\optfull} (from 67.0\% to 59.0\%), where the ability to recall historical code fixes is vital for solving real-world instances.
Finally, the \textbf{Debate Mechanism}, together with the \textbf{Debate Memory}, consistently improves performance over the single-agent team, validating that decentralized cross-verification effectively mitigates single-model myopia problems.

\subsection{In-Depth Analysis on Debate}

\subsubsection{Decentralized Debate vs. Centralized Judge}


%

To evaluate the effectiveness of decentralized debate against selection methods based on a third-party judge, we conduct a comparative study of two aggregation paradigms for combining multiple candidate solutions. The candidates are generated by two separate instances of the initial solution-generation pipeline used in \texttt{Agora-Opt} with solution memory, namely the \textbf{GPT-4o Team} and the \textbf{DeepSeek-V3 Team}. Built on these same two teams, we compare three aggregation strategies. In our decentralized debate framework, the two teams interact to reach a final solution. In \textbf{GenPRM-Selection} and \textbf{DeepSeek-Selection}, this interaction is replaced by a centralized judge that evaluates the two completed candidate solutions and selects one as the final answer, where the judges are instantiated by \texttt{DeepSeek-V3} and \texttt{GenPRM}, respectively, and \texttt{GenPRM} is the state-of-the-art generative process reward model for LLM-based optimization modeling proposed by \citet{zhou2025steporlm}. 
Therefore, all three methods share the same candidate-generation pipeline and differ only in the aggregation mechanism.

\begin{figure*}[ht]
  \centering
  \includegraphics[width=0.8\linewidth]{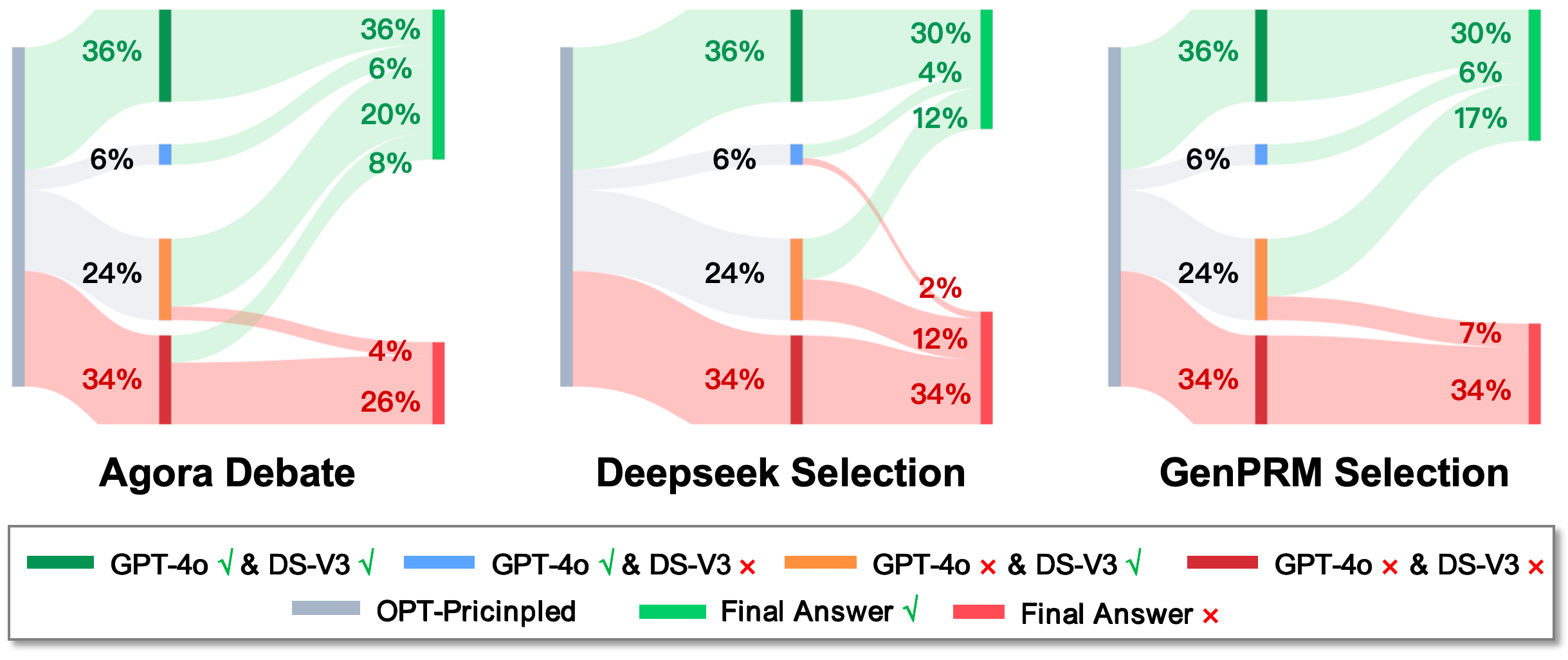}
  \caption{Sankey diagrams of outcome transitions on \texttt{\optfull} for Agora Debate, DeepSeek Selection, and GenPRM Selection. The flows group instances by the initial correctness pattern of the \texttt{GPT-4o} Team and the \texttt{DeepSeek-V3} Team, and trace their transitions to final correct or incorrect answers under each method.}
  \label{fig:sankey}
\end{figure*}

We evaluate these methods on three challenging benchmarks, \texttt{ComplexLP}, \texttt{IndustryOR}, and \texttt{\optfull}, as reported in \cref{table:debate_comparison}. To further examine how the aggregation mechanism affects the transition from initial team outputs to final correctness, we also visualize the instance-level state transitions on \texttt{\optfull} using the Sankey diagram in \cref{fig:sankey}. Specifically, we partition instances into four categories according to the correctness of the two initial teams: both correct, only \texttt{GPT-4o} correct, only \texttt{DeepSeek-V3} correct, or both wrong. We then trace how each aggregation strategy maps these four starting conditions to a correct or incorrect final answer.

\input{tables/debate_beyond_selection}

Our results reveal three clear advantages of decentralized debate over centralized judge selection.

\textbf{Breaking the Selection Bottleneck.}
Centralized judges are fundamentally limited to selecting among the candidate answers already produced by the two teams. Their performance therefore depends not only on the quality of the candidate pool, but also on the judge’s ability to identify the better solution. As shown in \cref{table:debate_comparison}, on \texttt{\optfull}, the stronger individual team \texttt{DeepSeek-V3} achieves 63.0\% accuracy, while \texttt{DeepSeek-Selection} and \texttt{GenPRM-Selection} reach only 52.0\% and 59.0\%, respectively. A similar pattern appears on \texttt{IndustryOR}, where \texttt{DeepSeek-V3} attains 71.0\%, compared with 65.0\% for \texttt{DeepSeek-Selection} and 70.0\% for \texttt{GenPRM-Selection}. By contrast, Agora-Opt consistently outperforms both the centralized judges and the best individual team across all three benchmarks, reaching 81.1\% on \texttt{ComplexLP}, 74.0\% on \texttt{IndustryOR}, and 70.0\% on \texttt{\optfull}. This suggests that decentralized debate does more than select between existing candidates: it enables the system to refine and improve upon the initial candidate pool itself.

\textbf{Debate Enables Self-Correction Beyond Static Selection.}
The most important structural difference appears in the hardest cases, where neither initial team is correct. As shown in \cref{fig:sankey}, on \texttt{\optfull}, 34\% of instances fall into the \emph{both-wrong} bucket. For centralized selection methods, these cases are effectively unrecoverable: when both candidate solutions are wrong, a judge can only select one wrong answer. 
By contrast, decentralized debate in \texttt{Agora-Opt} uniquely recovers 8\% of all instances from the both-wrong bucket to final correct answers, that is, 8/34 (about 23.5\%) of the dual-failure cases. 
No such recovery appears for either \texttt{DeepSeek-Selection} or \texttt{GenPRM-Selection}.
This result shows that decentralized debate is generative rather than merely selective: through active debate, the two teams can synthesize partial insights, identify mutual logical flaws, and collaboratively produce a correct solution that neither could achieve independently.

\textbf{Reliable Use of Complementary Strengths.}
Decentralized debate is also more effective when the two teams disagree, that is, when one initial solution is correct and the other is not. On \texttt{\optfull}, these asymmetric cases account for 30\% of all instances in total (6\% where only \texttt{GPT-4o} is correct and 24\% where only \texttt{DeepSeek-V3} is correct). In these cases, \texttt{Agora-Opt} converts 26/30 of them into final correct answers, compared with 23/30 for \texttt{GenPRM-Selection} and only 16/30 for \texttt{DeepSeek-Selection}. This suggests that, rather than making a one-shot choice between two static outputs, decentralized debate more reliably preserves the correct reasoning when only one team is correct and allows the two teams to leverage complementary partial strengths through critique and revision. 

Taken together, these results show that the core advantage of Agora-Opt lies in turning disagreement into productive interaction, producing a more robust final solution than centralized judge selection.

\subsubsection{Impact of Debate Rounds}

\begin{figure*}[htpb]
  \centering
  \includegraphics[width=1\linewidth]{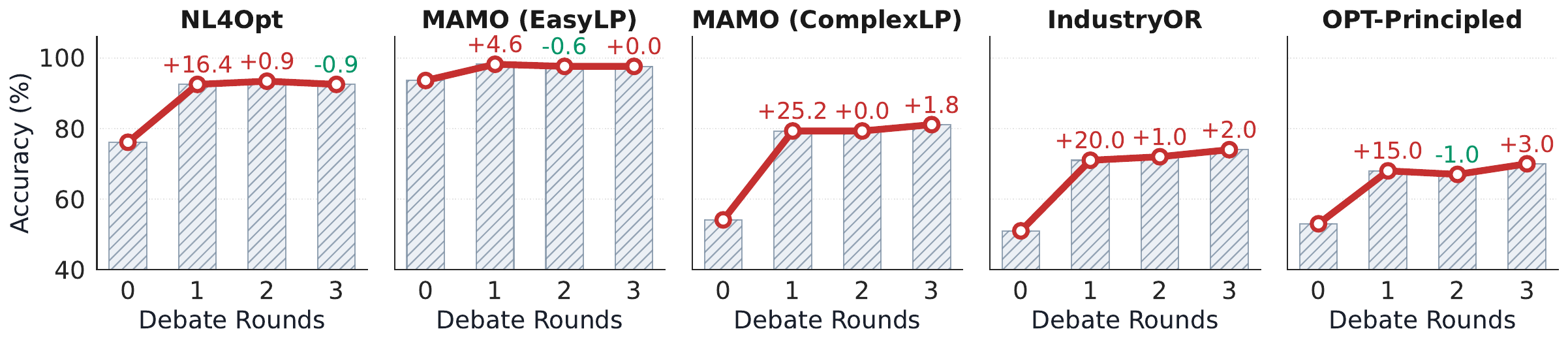}
  \caption{The performance improvement w.r.t. the number of debate rounds.}
  \label{fig:debate round}
\end{figure*}

We investigate the impact of the number of debate rounds on the performance. 
\cref{fig:debate round} illustrates the accuracy trajectory across four benchmarks as the number of debate rounds increases from 0 to 3.
Two key patterns emerge from these results regarding the efficiency and necessity of the debate mechanism.

\textbf{Rapid Convergence and Efficiency.}
The most significant performance leap occurs consistently in the first round of debate ($T_\text{max}=0 \rightarrow 1$).
As shown in \cref{fig:debate round}, Agora-Opt achieves the vast majority of its improvement immediately after the first round of debate.
For instance, on \texttt{MAMO(ComplexLP)}, the accuracy surges by 25.2\% (from 54.1\% to 79.3\%) in a single round.
Subsequent rounds ($T_\text{max}=2, 3$) yield marginal gains or stabilize the performance (e.g., \texttt{IndustryOR} rises from 71.0\% to 74.0\%, while \texttt{NL4Opt} fluctuates within 1\%).
This ``one-round convergence'' property suggests that our framework is computationally efficient, as the majority of formulation errors and constraint misses are identified and corrected during the initial critique, rendering prolonged debates unnecessary for most instances.

\textbf{The Risk of ``Over-thinking'' on Simple Tasks.}
While extended debate consistently benefits complex tasks (e.g., \texttt{IndustryOR} improves monotonically up to round 3), it can be counterproductive on simpler benchmarks.
As observed in \texttt{NL4Opt}, \texttt{MAMO(EasyLP)} and \texttt{\optfull}, performance fluctuates or even regresses in later rounds (e.g., \texttt{NL4Opt} drops from 93.4\% at $T_\text{max}=2$ to 92.5\% at $T_\text{max}=3$).
This phenomenon aligns with recent observations in reasoning models where ``over-thinking'' simple problems can introduce unnecessary noise or hallucinatory corrections into an already correct solution~\citep{chennot,liu2025think}.
This finding suggests that a fixed number of debate rounds may not be optimal for all problem types.
Consequently, developing an adaptive mechanism to dynamically control the debate depth based on instance complexity represents a promising direction for future work.

\subsubsection{Generalization of Decentralized Debate Procotol}
\label{subsec:swap_methods}

\input{tables/debate_generalization}

We apply the decentralized debate protocol beyond the Agora-Opt framework to other representative baselines to assess its portability and generalization. 
\cref{tab:swap methods} reports the performance of applying our debate mechanism across two distinct paradigms: agentic methods (e.g., OptiMUS-v0.3~\citep{ahmaditeshnizi2024optimus}, \texttt{CAFA}~\citep{deng24cafa}) and training-centric methods (e.g., ORLM~\citep{huang2025orlm}, StepORLM~\citep{zhou2025steporlm}). 
The results illustrate three key findings regarding the generalization capability and the performance dynamics of our decentralized debate protocol.

\textbf{Debate as Method-agnostic Principle.}
Across both paradigms, injecting decentralized debate improves performance without modifying the underlying solver architecture or training pipeline, confirming that agentic debate functions as a high-level meta-strategy effective across diverse architectures. Concretely, in the closed-source setting, enabling debate increases the macro-average from 58.9\% (\texttt{OptiMUS}) and 47.9\% (\texttt{CAFA}) to 70.7\%. Similarly, in the open-source 8B setting, by layering our debate protocol onto the static, training-centric StepORLM, we can also observe a performance boost (from 71.4\% to 76.3\%). These results demonstrate that our proposed debate mechanism is not merely an algorithmic component of Agora-Opt but a modular plug-in layer that effectively generalizes to external methods, unlocking their collective intelligence for optimization modeling.


\textbf{Performance Gains on Hard Tasks vs. Simple Tasks.}
A consistent pattern is that debate tends to deliver the largest benefits on benchmarks where modeling mistakes and verification failures are more likely, while offering limited gains, or even mild regressions, on benchmarks that are already close to saturation. In the open-source 8B setting, debate substantially improves the more challenging benchmarks, including \texttt{CompOR} (55.6\%$\rightarrow$66.7\%), \texttt{IndOR} (40.0\%$\rightarrow$53.0\%), and \texttt{\optfull} (21.0\%$\rightarrow$32.0\%), while slightly reducing performance on easier datasets such as \texttt{NL4OPT} (97.7\%$\rightarrow$96.7\%) and \texttt{MAMO(EasyLP)} (97.2\%$\rightarrow$96.7\%). This behavior also aligns with our earlier debate-round sensitivity analysis: debate is most valuable when it helps surface and correct long-tail formulation errors, but additional critique can introduce unnecessary revisions when the base solution is already highly reliable.

\textbf{Diminishing Marginal Utility.} 
A comparative analysis reveals a ``diminishing returns'' trend relative to the base method's capability. For the agentic methods, wrapping \texttt{OptiMUS} with debate yields a large macro-average gain of +11.8\%. In contrast, when applied to the stronger training-centric baseline \texttt{StepORLM}, debate still improves performance but with a smaller gain of +4.9\%. This suggests that while decentralized debate is highly effective at correcting reasoning errors in weaker or intermediate models, its marginal impact naturally stabilizes as the base system approaches the solvability ceiling of current benchmarks.



%% file: tables/LLM_compatibility.tex
\begin{table*}[htpb]
\centering
\caption{Performance comparison of Agora-Opt with different backbone LLM combinations.}
\label{table:agora_variants}

\resizebox{\textwidth}{!}{%
\begin{tabular}{l ccccccccc}
\toprule
Model & NL4OPT & \multicolumn{2}{c}{MAMO} & NLP4LP & CompOR & IndOR & ReSoc. & \optshort & Avg. \\
\cmidrule(lr){3-4}
 & & Easy & Complex & & & & & & \\
\midrule

\texttt{Gemini-2.5-Pro}
& 82.6 & 87.9 & 52.3 & 94.9 & 66.7 & 63.0 & 89.6 & 57.0 & 74.2 \\
\texttt{DeepSeek-V3} 
& 79.8 & 95.2 & 53.2 & 92.1 & 55.6 & 56.0 & 85.1 & 26.0 & 67.9 \\
\texttt{GPT-4o}
& 76.1 & 93.6 & 54.1 & 91.6 & 44.4 & 51.0 & 84.1 & 16.0 & 63.9 \\

\midrule

\texttt{Agora-Opt (GPT-4o + DS-V3)} 
& 92.5 & 97.6 & 81.1 & 98.9 & 66.7 & 74.0 & 96.3 & 70.0 & 84.6 \\

\texttt{Agora-Opt (DS-V3 + Gemini)}  
& 92.5 & 96.0 & 82.0 & 98.9 & 66.7 & 76.0 & 95.8 & 75.0 & 85.4\\

\texttt{Agora-Opt (GPT-4o + Gemini)}  
& 92.5 & 96.9 & 84.7 & 99.4 & 66.7 & 83.0 & 95.8 & 63.0 & 85.3\\

\bottomrule
\end{tabular}%
}

\par
\scriptsize
\vspace{2pt}
\textbf{Note}: DS-V3: DeepSeek-V3, Gemini: Gemini-2.5-Pro. IndOR: IndustryOR. ReSoc.: ReSocratic, \optshort: \optfull.

\end{table*}









%% file: tables/ablation_study.tex
\begin{table}[htpb]
\centering
\caption{Ablation study of different variants on ComplexLP, IndustryOR and \optfull\ benchmarks. Starting with the full Agora-Opt, we sequentially remove key components, gradually reducing the model to a vanilla LLM.}
\label{table:ablation_study}

\begin{tabular}{lccc}
\toprule
\textbf{Variants} & \textbf{ComplexLP} & \textbf{IndustryOR} & \textbf{\optfull}\\
\midrule

Agora-Opt & 81.1 & 74.0 & 70.0\\ 
\midrule

Remove Debate Memory & 80.2 & 72.0 & 67.0\\
Remove Debug Memory & 76.6 & 64.0 & 59.0\\
Remove Solution Memory & 64.9 & 61.0 & 56.0 \\
Remove Debate (Single Agent Team) & 62.2 & 59.0 & 21.0\\
Remove Agent Team (Vanilla GPT-4o) & 54.1 & 51.0 & 16.0\\

\bottomrule
\end{tabular}
\end{table}

%% file: tables/debate_beyond_selection.tex
\begin{table}[htpb]
\centering
\caption{Comparsion Study of Debate mechanism vs. GenRPM selection and Deepseek selection. We conducted the experiment on ComplexLP, IndustryOR and \optfull\ benchmarks. }
\label{table:debate_comparison}

\begin{tabular}{lccc}
\toprule
\textbf{Model} & \textbf{ComplexLP} & \textbf{IndustryOR} & \textbf{\optfull}\\

\midrule
GPT-4o Team & 79.3 & 60.0 & 32.0 \\
Deepseek-V3 Team & 78.4 & 71.0 & 63.0 \\
\midrule

Agora-Opt & 81.1 & 74.0 & 70.0\\ 
Deepseek-Selection & 80.2 & 65.0 & 52.0 \\
GenPRM-Selection & 79.3 & 70.0 & 59.0\\ 

\bottomrule

\end{tabular}

\par
\scriptsize
\vspace{2pt}
\textbf{Note}: Team here means our initial solution generation method used in Agora-Opt, working with the solution memory.

\end{table}

%% file: tables/debate_generalization.tex
\begin{table*}[htpb]
\footnotesize
\centering
\caption{Performance comparison of varying models under debate settings. 
The table is divided into two groups: Closed-source model baselines and Open-source (8B) model baselines.
Best results are highlighted in \textbf{bold} and second-highest values are \underline{underlined}.}
\label{tab:swap methods}

\setlength{\tabcolsep}{2pt} 
\adjustbox{max width=1.0\textwidth}{
\begin{tabularx}{\textwidth}{l c *{9}{>{\centering\arraybackslash}X}}
\toprule
Model & Params & NL4OPT & \multicolumn{2}{c}{MAMO} & NLP4LP & CompOR & IndOR & ReSoc. & \optshort & Avg. \\
\cmidrule(lr){4-5}
 & & & EasyLP & ComplexLP & & & & & \\
\midrule

\rowcolor{gray_header}
\multicolumn{11}{c}{\textbf{\textit{Debate on Agentic Methods}}} \\
\addlinespace[2pt]

\texttt{OptiMUS(GPT-4o)} & Closed
& \underline{76.2}
& \underline{78.0}
& \underline{46.8}
& \underline{88.8}
& \underline{44.1}
& \underline{30.0}
& \underline{87.6}
& 20.0
& \underline{58.9}
\\

\texttt{CAFA(GPT-4o)}
& Closed
& 68.1
& 71.2
& 44.5
& 50.0
& 46.4
& 39.0
& 40.1
& \underline{24.0}
& 47.9
\\

\rowcolor{purple_highlight}
\texttt{Debate} & Closed
& \textbf{85.5} & \textbf{94.5} & \textbf{65.8} & \textbf{92.7} & \textbf{55.6} & \textbf{46.0} & \textbf{90.1} & \textbf{35.0} & \textbf{70.7} \\

\midrule
\midrule

\rowcolor{gray_header}
\multicolumn{11}{c}{\textbf{\textit{Debate on Training-centric Methods}}} \\
\addlinespace[2pt]

\texttt{ORLM} & 8B
& 73.8
& 90.4
& 59.5
& 76.4
& 50.0
& 38.0
& 61.8
& 19.0
& 58.6 \\

\texttt{StepORLM} & 8B
& \textbf{97.7}
& \textbf{97.2}
& \textbf{79.3}
& \underline{97.8}
& \underline{55.6}
& \underline{40.0}
& \underline{82.6}
& \underline{21.0}
& \underline{71.4}
\\

\rowcolor{purple_highlight}
\texttt{Debate} & 8B 
& \underline{96.7} & \underline{96.7} & \underline{78.4} & \textbf{98.3} & \textbf{66.7} & \textbf{53.0} & \textbf{88.2} & \textbf{32.0} & \textbf{76.3} \\

\bottomrule
\end{tabularx}
}
\par
\vspace{2pt}
\scriptsize
\textbf{Abbreviations}: OptiMUS: OptiMUS-v0.3, CompOR: ComplexOR, IndOR: IndustryOR, ReSoc.: ReSocratic, \optshort: \optfull, Avg: Macro-Average.
\end{table*}

%% file: appendix/benchmarks.tex
\section{Public Benchmarks}\label{appendix:benchmarks}

A diverse array of datasets and construction methodologies has been proposed to enhance the capability of Large Language Models (LLMs) in addressing Operations Research (OR) problems, with these benchmarks becoming increasingly sophisticated and realistic \citep{yang2025optibenchmeetsresocraticmeasure, AhmadiTeshnizi2024, yang2024optibenchmeetsresocraticmeasure}. A recent survey by \citet{xiao2025survey} provides a comprehensive and timely review of benchmark datasets for optimization problems described in natural language. This study identified prevalent issues within existing datasets, including logical fallacies, ambiguous definitions, and incorrect ground-truth answers, and rectified a portion of these errors. Building upon these initial efforts, our work involves further manual verification and cleaning of these widely-used benchmarks. Subsequently, the data is processed into a unified format, establishing a more robust and standardized foundation for subsequent research and model assessment.

\subsection{NL4Opt}
\label{sec:nl4opt}
The NL4Opt dataset was introduced as part of the NeurIPS 2022 NL4Opt competition \citep{ramamonjison2022augmenting, Ramamonjison2023}. It comprises approximately 1,100 annotated word problems describing linear programming (LP) scenarios, primarily aimed at bridging the gap between natural language descriptions and formal optimization models. \citet{xiao2025survey} reported that NL4Opt possesses a complexity score of 5.59 and an original error rate of at least 26.4\%. Each problem presents a self-contained, real-world scenario (e.g., resource allocation, scheduling) with annotations identifying key optimization elements. The official shared task is divided into two subtasks: (1) recognizing and labeling problem entities, and (2) generating the corresponding mathematical model. While the original dataset provides 713 training, 99 validation, and 289 test instances, we selected 213 validated problems from this benchmark for our evaluation.

\subsection{MAMO}
\label{sec:mamo}
The MAMO dataset was introduced by \citet{huang2024mamo} as a benchmark for evaluating LLMs on complex mathematical modeling tasks, focusing specifically on the formulation process rather than mere solution accuracy. The dataset concentrates on linear programming (LP) and mixed-integer linear programming (MILP) problems, excluding nonlinear or differential equation modeling. MAMO is categorized into two main components: Easy LP (652 instances) and Complex LP (211 instances). 
Analysis by \citet{xiao2025survey} assigned EasyLP a complexity of 7.12 with an error rate of at least 8.13\%, while ComplexLP exhibits a significantly higher complexity of 13.35 and an error rate of at least 23.7\%. By delegating the solving step to an external OR solver, MAMO isolates whether an LLM can accurately translate detailed natural language descriptions into correct mathematical formulations. From this benchmark, we selected 545 easy and 111 complex instances for our study.

\subsection{NLP4LP}
\label{sec:nlp4lp}
The NLP4LP benchmark, introduced alongside the OptiMUS system~\citep{ahmaditeshnizi2024optimus}, evaluates LLM-based agents on translating natural language OR problems into solver-ready code and mathematical models. The survey by \citet{xiao2025survey} assigned it a complexity score of 5.58 and identified an error rate of at least 21.7\% in the original version. It consists of 269 human-authored LP and MILP problems, each representing an optimization scenario from classical OR domains such as scheduling, knapsack allocation, and production planning. This benchmark serves as a fine-grained testbed for assessing formulation accuracy. In our work, we chose to evaluate 178 verified instances from NLP4LP.

\subsection{ComplexOR}
\label{sec:complexor}
The ComplexOR dataset was introduced  as part of the Chain-of-Experts framework to stress-test the reasoning and modeling capabilities of LLMs on scenarios more demanding than those in earlier datasets \citep{xiao2024chainofexperts}. The analysis by \citet{xiao2025survey} noted a complexity score of 5.98 and a substantial error rate of at least 24.3\%. ComplexOR was curated from advanced OR case studies and challenging optimization problem classes, such as multi-step lot-sizing with setup costs, intricate scheduling, and supply chain design. These problems are primarily large-scale or conceptually complex MILPs requiring multi-step reasoning and often involving hierarchical or time-indexed structures. The dataset is designed to facilitate the development and evaluation of advanced reasoning strategies for OR problem-solving. Following our verification process, we selected 18 reliable problems from this dataset.

\subsection{IndustryOR}
\label{sec:industryor}
The IndustryOR dataset, introduced by \citet{huang2025orlm}, is the first industrial-scale operations research benchmark for LLMs. It comprises 100 real-world optimization scenarios spanning sectors such as manufacturing, logistics, finance, and energy. The study by \citet{xiao2025survey} highlighted the high difficulty of this dataset, reporting a complexity score of 14.06 and an error rate of at least 54.0\%. Designed to evaluate LLMs on practical, domain-specific tasks, IndustryOR covers five OR categories, including linear, integer, mixed-integer, nonlinear, and others, categorized by difficulty (Easy, Medium, Hard). For our experiments, we selected 42 validated problems from this benchmark.

\subsection{ReSocratic}
\label{sec:resocratic}
The ReSocratic dataset was introduced by \citet{yang2025optibenchmeetsresocraticmeasure}, representing a data synthesis method proposed in the same work. This method synthesizes formatted optimization demonstrations in a reverse manner before back-translating them into natural language questions. By leveraging these intermediate reasoning steps, ReSocratic achieves higher quality than previous synthesis methods. In this paper, we selected 403 validated problems from the original ReSocratic-29K dataset.

%% file: appendix/prompts.tex
\section{Prompts for Agora-Opt}
\label{appendix_prompts}

In this section, we provide the detailed prompts used in our Agora-opt framework.

\subsection{Formulator Agent}

\begin{promptbox}{Formulator Prompt}

\textbf{Role \& Goal}\\
You are the Formulator in an optimization agent team. Your job is to turn a natural-language OR problem into a clear, solver-ready mathematical formulation.

\textbf{Deliverable (strict)}\\
Output ONLY a \code{<formulation>...</formulation>} block containing:
\begin{itemize}
    \item Sets / indices
    \item Parameters (with meanings/units if relevant)
    \item Decision variables (with domains/types)
    \item Objective function (mark min or max explicitly)
    \item Constraints (well-typed, no prose-only statements)
\end{itemize}

\textbf{Requirements \& Guardrails}
\begin{itemize}
    \item Be explicit and unambiguous; avoid loose prose.
    \item Name every symbol you introduce; avoid dangling references.
    \item If data hints suggest bounds, include them.
    \item Keep solver-ready structure (no storytelling or extra commentary).
\end{itemize}

\textbf{Memory Context}\\
The block is retrieved similar solved cases or experiences, which are only used as references for this task.\\
\code{\{\{solution\_memory\}\}}

\textbf{Problem Statement}\\
\code{\{\{problem\}\}}
\end{promptbox}

\subsection{Programmer Agent}

\begin{promptbox}{Programmer Prompt}
\textbf{Role \& Goal}\\
You are the Programmer in an optimization agent team. Translate the given formulation into executable Gurobi Python code.

\textbf{Deliverable (strict)}\\
Output ONLY a \code{<python>...</python>} block with the FULL runnable script.

\textbf{Guardrails}
\begin{itemize}
    \item Use gurobipy only (no pulp/cvxpy/coptpy).
    \item Define every symbol you reference (sets, parameters, big-M, helper dicts).
    \item If big-M is used, pick explicit, defensible bounds.
    \item Ensure the script is self-contained (imports, data parsing, model solve).
    \item Print the objective value when solved via: \code{print("OBJECTIVE\_VALUE:", model.objVal)}.
    \item Avoid placeholder code or ellipses; produce complete, runnable Python.
\end{itemize}

\textbf{Memory Context}\\
The block is retrieved similar solved cases or experiences, which are only used as references for this task.\\
\code{\{\{solution\_memory\}\}}

\textbf{Inputs}\\
Problem: \code{\{\{problem\}\}}\\
Formulation: \code{\{\{formulation\}\}}
\end{promptbox}

\subsection{Debugger Agent}

\begin{promptbox}{Debugger Prompt}
\textbf{Role \& Task}\\
You are the Debugger in an optimization agent team. The generated code failed to run—diagnose the root cause and return a fixed, runnable script.

\textbf{Deliverable (strict)}\\
Output ONLY a \code{<python>...</python>} block with the FULL corrected script.

\textbf{Requirements \& Guardrails}
\begin{itemize}
    \item Preserve the formulation intent; change only what is necessary to make it run correctly.
    \item Ensure all symbols are defined before use; no missing imports or variables.
    \item If big-M is used, choose explicit, defensible bounds.
    \item Make sure the code prints the objective value when solved.
\end{itemize}

\textbf{Inputs}\\
Problem: \code{\{\{problem\}\}}\\
Formulation: \code{\{\{formulation\}\}}\\
Failed code: \code{\{\{current\_code\}\}}\\
Execution status: \code{\{\{exec\_status\}\}}\\
Solver flag: \code{\{\{solver\_flag\}\}}\\
Stderr: \code{\{\{stderr\}\}}

\textbf{Historical Debug Memory}\\
\code{\{\{debug\_memory\}\}}
\end{promptbox}

\subsection{Retrieval Analysis}

\begin{promptbox}{Retrieval Analysis Prompt}
You are an expert in optimization modeling. You will analyze multiple similar solved problems to extract \textbf{transferable insights} for a new problem.

\textbf{Current Problem to Solve:}\\
\code{\{current\_problem\_desc\}}

\textbf{Retrieved Similar Cases (Complete):}\\
\code{\{full\_cases\}}

\textbf{Your Task:}\\
Analyze ALL the cases above \textbf{holistically} and provide a structured analysis that will guide solving the current problem.

\textbf{Focus on:}
\begin{enumerate}
    \item \textbf{Problem Type \& Structure}: What category do these problems fall into? (e.g., production planning, resource allocation, scheduling, network flow)
    \item \textbf{Common Modeling Patterns}: 
    \begin{itemize}
        \item What decision variables are typically used?
        \item What types of constraints appear repeatedly?
        \item How are objectives typically formulated?
    \end{itemize}
    \item \textbf{Key Techniques \& Tricks}:
    \begin{itemize}
        \item Any specific Gurobi features? (e.g., \code{addConstrs}, \code{quicksum}, binary variables, \code{setParam})
        \item Modeling tricks? (e.g., big-M, indicator constraints, piecewise linear)
        \item Data structure patterns? (e.g., dictionaries for indices, list comprehensions)
    \end{itemize}
    \item \textbf{Adaptation Guidance}:
    \begin{itemize}
        \item What aspects of the current problem are similar to the retrieved cases?
        \item What's different and requires new thinking?
        \item Which parts of the solution approaches can be directly adapted?
    \end{itemize}
\end{enumerate}

\textbf{Output Format}: \\
Provide a concise, actionable analysis (300-500 words) structured by the 4 points above. Be specific with code patterns and techniques, not just high-level descriptions.
\end{promptbox}

\subsection{Debate Agent}

\begin{promptbox}{Debate Prompt}
\textbf{Role \& Goal}\\
You are participating in an agentic debate to improve an optimization model. You see the problem, your candidate, and the other team's candidate. Critique both and return a revised formulation + code for your own candidate.

\textbf{Deliverable (strict)}\\
Output ONLY the following blocks, in order:\\
\code{<formulation>...</formulation>}\\
\code{<python>...</python>}

\textbf{Guardrails}
\begin{itemize}
    \item Keep the objective/constraints logically consistent with the problem.
    \item Incorporate useful ideas from the other candidate when they improve correctness.
    \item Fix execution/stability issues you spot (missing symbols, bounds, solver flags).
    \item If big-M is needed, choose explicit, defensible bounds.
    \item Ensure the code is self-contained and runnable with gurobipy.
    \item Preserve clarity: no placeholder text, no truncated code.
\end{itemize}

\textbf{Inputs}\\
Problem: \code{\{\{problem\}\}}\\
Your candidate: Formulation: \code{\{\{my\_formulation\}\}}, Code: \code{\{\{my\_code\}\}}\\
Other candidate: Formulation: \code{\{\{other\_formulation\}\}}, Code: \code{\{\{other\_code\}\}}

\textbf{Debate Memory}\\
Debate memory contains retrieved summaries of past debates where disagreements were resolved. Use these to guide how you reconcile differing candidates: note decisive arguments, guardrails, and modeling patterns; adapt, don't copy blindly:\\
\code{\{\{debate\_memory\}\}}
\end{promptbox}

\subsection{Debate Memory Summarizer}

\begin{promptbox}{Debate Memory Summarization Prompt}
You are helping an optimisation-debate memory builder.

\textbf{Inputs}\\
Problem description: \code{\{description\}}\\
Initial disagreement: Model A result: \code{\{initial\_A\_result\}}, Model B result: \code{\{initial\_B\_result\}}\\
Debate transcript: \code{\{history\_text\}}\\
Final consensus objective: \code{\{final\_result\}}

Please return a JSON object with the following keys:
\begin{itemize}
    \item "summary": 2-3 sentences explaining how the debate resolved the mismatch.
    \item "mismatch\_reason": concise reason for the disagreement.
    \item "decisive\_argument": specific insight that convinced both sides.
    \item "guardrails": list of actionable bullet points the next debater should follow.
    \item "modeling\_patterns": list of reusable modeling tricks/structures that appeared.
\end{itemize}

JSON ONLY. No prose outside the JSON.
\end{promptbox}

%% file: custom.bib
@article{ramamonjison2022augmenting,
  title={Augmenting operations research with auto-formulation of optimization models from problem descriptions},
  author={Ramamonjison, Rindranirina and Li, Haley and Yu, Timothy T and He, Shiqi and Rengan, Vishnu and Banitalebi-Dehkordi, Amin and Zhou, Zirui and Zhang, Yong},
  journal={arXiv preprint arXiv:2209.15565},
  year={2022}
}

@inproceedings{xiao2023chain,
  title={Chain-of-Experts: When LLMs Meet Complex Operations Research Problems},
  author={Xiao, Ziyang and Zhang, Dongxiang and Wu, Yangjun and Xu, Lilin and Wang, Yuan Jessica and Han, Xiongwei and Fu, Xiaojin and Zhong, Tao and Zeng, Jia and Song, Mingli and others},
  booktitle={The Twelfth International Conference on Learning Representations},
  year={2023}
}

@inproceedings{nl4opt,
  title={{NL4Opt} competition: Formulating optimization problems based on their natural language descriptions},
  author={Ramamonjison, Rindranirina and Yu, Timothy and Li, Raymond and Li, Haley and Carenini, Giuseppe and Ghaddar, Bissan and He, Shiqi and Mostajabdaveh, Mahdi and Banitalebi-Dehkordi, Amin and Zhou, Zirui and others},
  booktitle={NeurIPS 2022 Competition Track},
  pages={189--203},
  year={2023},
  organization={PMLR}
}

@misc{ahmaditeshnizi2024optimus,
      title={OptiMUS-0.3: Using Large Language Models to Model and Solve Optimization Problems at Scale}, 
      author={Ali AhmadiTeshnizi and Wenzhi Gao and Herman Brunborg and Shayan Talaei and Madeleine Udell},
      year={2024},
      eprint={2407.19633},
      archivePrefix={arXiv},
      primaryClass={cs.AI},
      url={https://arxiv.org/abs/2407.19633}, 
}

@misc{yang2024optibenchmeetsresocraticmeasure,
      title={OptiBench Meets ReSocratic: Measure and Improve LLMs for Optimization Modeling}, 
      author={Zhicheng Yang and Yiwei Wang and Yinya Huang and Zhijiang Guo and Wei Shi and Xiongwei Han and Liang Feng and Linqi Song and Xiaodan Liang and Jing Tang},
      year={2024},
      eprint={2407.09887},
      archivePrefix={arXiv},
      primaryClass={cs.LG},
      url={https://arxiv.org/abs/2407.09887}, 
}

@inproceedings{AhmadiTeshnizi2024,
  author = {Ali AhmadiTeshnizi and Wenzhi Gao and Madeleine Udell},
  title = {{OptiMUS: Scalable Optimization Modeling with (MI)LP Solvers and Large Language Models}},
  booktitle = {Proceedings of the 41st International Conference on Machine Learning (ICML)},
  year = {2024}
}

@article{Ramamonjison2023,
  author = {Rindranirina Ramamonjison and Timothy T. Yu and Raymond Li and Haley Li and Giuseppe Carenini and Bissan Ghaddar and Shiqi He and Mahdi Mostajabdaveh and Amin Banitalebi-Dehkordi and Zirui Zhou and Yong Zhang},
  title = {{NL4Opt Competition: Formulating Optimization Problems Based on Their Natural Language Descriptions}},
  journal = {arXiv preprint arXiv:2303.08233},
  year = {2023}
}

@article{chen2025stepwise,
  title={Stepwise Guided Policy Optimization: Coloring your Incorrect Reasoning in GRPO},
  author={Chen, Peter and Li, Xiaopeng and Li, Ziniu and ChenD, Xi and Lin, Tianyi},
  journal={arXiv preprint arXiv},
  volume={2505},
  year={2025}
}

@article{xiao2025survey,
  title={A Survey of Optimization Modeling Meets LLMs: Progress and Future Directions},
  author={Xiao, Ziyang and Xie, Jingrong and Xu, Lilin and Guan, Shisi and Zhu, Jingyan and Han, Xiongwei and Fu, Xiaojin and Yu, WingYin and Wu, Han and Shi, Wei and others},
  journal={arXiv preprint arXiv:2508.10047},
  year={2025}
}

@misc{yang2025optibenchmeetsresocraticmeasure,
      title={OptiBench Meets ReSocratic: Measure and Improve LLMs for Optimization Modeling}, 
      author={Zhicheng Yang and Yiwei Wang and Yinya Huang and Zhijiang Guo and Wei Shi and Xiongwei Han and Liang Feng and Linqi Song and Xiaodan Liang and Jing Tang},
      year={2025},
      eprint={2407.09887},
      archivePrefix={arXiv},
      primaryClass={cs.LG},
      url={https://arxiv.org/abs/2407.09887}, 
}

@misc{chen2025solverinformedrlgroundinglarge,
      title={Solver-Informed RL: Grounding Large Language Models for Authentic Optimization Modeling}, 
      author={Yitian Chen and Jingfan Xia and Siyu Shao and Dongdong Ge and Yinyu Ye},
      year={2025},
      eprint={2505.11792},
      archivePrefix={arXiv},
      primaryClass={cs.AI},
      url={https://arxiv.org/abs/2505.11792}, 
}

@misc{yang2025qwen3technicalreport,
      title={Qwen3 Technical Report}, 
      author={An Yang and Anfeng Li and Baosong Yang and Beichen Zhang and Binyuan Hui and Bo Zheng and Bowen Yu and Chang Gao and Chengen Huang and Chenxu Lv and Chujie Zheng and Dayiheng Liu and Fan Zhou and Fei Huang and Feng Hu and Hao Ge and Haoran Wei and Huan Lin and Jialong Tang and Jian Yang and Jianhong Tu and Jianwei Zhang and Jianxin Yang and Jiaxi Yang and Jing Zhou and Jingren Zhou and Junyang Lin and Kai Dang and Keqin Bao and Kexin Yang and Le Yu and Lianghao Deng and Mei Li and Mingfeng Xue and Mingze Li and Pei Zhang and Peng Wang and Qin Zhu and Rui Men and Ruize Gao and Shixuan Liu and Shuang Luo and Tianhao Li and Tianyi Tang and Wenbiao Yin and Xingzhang Ren and Xinyu Wang and Xinyu Zhang and Xuancheng Ren and Yang Fan and Yang Su and Yichang Zhang and Yinger Zhang and Yu Wan and Yuqiong Liu and Zekun Wang and Zeyu Cui and Zhenru Zhang and Zhipeng Zhou and Zihan Qiu},
      year={2025},
      eprint={2505.09388},
      archivePrefix={arXiv},
      primaryClass={cs.CL},
      url={https://arxiv.org/abs/2505.09388}, 
}

@misc{qwen2025qwen25technicalreport,
      title={Qwen2.5 Technical Report}, 
      author={An Yang and Baosong Yang and Beichen Zhang and Binyuan Hui and Bo Zheng and Bowen Yu and Chengyuan Li and Dayiheng Liu and Fei Huang and Haoran Wei and Huan Lin and Jian Yang and Jianhong Tu and Jianwei Zhang and Jianxin Yang and Jiaxi Yang and Jingren Zhou and Junyang Lin and Kai Dang and Keming Lu and Keqin Bao and Kexin Yang and Le Yu and Mei Li and Mingfeng Xue and Pei Zhang and Qin Zhu and Rui Men and Runji Lin and Tianhao Li and Tianyi Tang and Tingyu Xia and Xingzhang Ren and Xuancheng Ren and Yang Fan and Yang Su and Yichang Zhang and Yu Wan and Yuqiong Liu and Zeyu Cui and Zhenru Zhang and Zihan Qiu},
      year={2025},
      eprint={2412.15115},
      archivePrefix={arXiv},
      primaryClass={cs.CL},
      url={https://arxiv.org/abs/2412.15115}, 
}

@inproceedings{wei2022chainofthought,
  title = {Chain-of-Thought Prompting Elicits Reasoning in Large Language Models},
  author = {Wei, Jason and Wang, Xuezhi and Schuurmans, Dale and Bosma, Maarten and ichter, brian and Xia, Fei and Chi, Ed and Le, Quoc V and Zhou, Denny},
  booktitle = {Advances in Neural Information Processing Systems},
  volume = {35},
  pages = {24824--24837},
  year = {2022},
  publisher = {Curran Associates, Inc.}
}

@techreport{openai_o3_2025,
  title={OpenAI o3 and o4-mini System Card},
  author={OpenAI},
  year={2025},
  institution={OpenAI},
  url={https://openai.com/index/o3-o4-mini-system-card/}
}

@article{liu2025think,
  title={Think How to Think: Mitigating Overthinking with Autonomous Difficulty Cognition in Large Reasoning Models},
  author={Liu, Yongjiang and Li, Haoxi and Ma, Xiaosong and Zhang, Jie and Guo, Song},
  journal={arXiv preprint arXiv:2507.02663},
  year={2025}
}

@inproceedings{chennot,
  title={Do NOT Think That Much for 2+ 3=? On the Overthinking of Long Reasoning Models},
  author={Chen, Xingyu and Xu, Jiahao and Liang, Tian and He, Zhiwei and Pang, Jianhui and Yu, Dian and Song, Linfeng and Liu, Qiuzhi and Zhou, Mengfei and Zhang, Zhuosheng and others},
  booktitle={Forty-second International Conference on Machine Learning}
}

@inproceedings{chen2025debatecoder,
  title={DebateCoder: Towards Collective Intelligence of LLMs via Test Case Driven LLM Debate for Code Generation},
  author={Chen, Jizheng and Du, Kounianhua and Dai, Xinyi and Zhang, Weiming and Wang, Xihuai and Wang, Yasheng and Tang, Ruiming and Zhang, Weinan and Yu, Yong},
  booktitle={Proceedings of the 63rd Annual Meeting of the Association for Computational Linguistics (Volume 1: Long Papers)},
  pages={12055--12065},
  year={2025}
}

@article{huang2025orlm,
  title={Orlm: A customizable framework in training large models for automated optimization modeling},
  author={Huang, Chenyu and Tang, Zhengyang and Hu, Shixi and Jiang, Ruoqing and Zheng, Xin and Ge, Dongdong and Wang, Benyou and Wang, Zizhuo},
  journal={Operations Research},
  year={2025},
  publisher={INFORMS}
}

@article{li2025swe,
  title={Swe-debate: Competitive multi-agent debate for software issue resolution},
  author={Li, Han and Shi, Yuling and Lin, Shaoxin and Gu, Xiaodong and Lian, Heng and Wang, Xin and Jia, Yantao and Huang, Tao and Wang, Qianxiang},
  journal={arXiv preprint arXiv:2507.23348},
  year={2025}
}

@article{long2024multi,
  title={Multi-expert prompting improves reliability, safety, and usefulness of large language models},
  author={Long, Do Xuan and Yen, Duong Ngoc and Luu, Anh Tuan and Kawaguchi, Kenji and Kan, Min-Yen and Chen, Nancy F},
  journal={arXiv preprint arXiv:2411.00492},
  year={2024}
}

@article{zhou2025memento,
  title={Memento: Fine-tuning llm agents without fine-tuning llms},
  author={Zhou, Huichi and Chen, Yihang and Guo, Siyuan and Yan, Xue and Lee, Kin Hei and Wang, Zihan and Lee, Ka Yiu and Zhang, Guchun and Shao, Kun and Yang, Linyi and others},
  journal={arXiv preprint arXiv:2508.16153},
  year={2025}
}

@article{irving2018ai,
  title={AI safety via debate},
  author={Irving, Geoffrey and Christiano, Paul and Amodei, Dario},
  journal={arXiv preprint arXiv:1805.00899},
  year={2018}
}

@inproceedings{du2023improving,
  title={Improving factuality and reasoning in language models through multiagent debate},
  author={Du, Yilun and Li, Shuang and Torralba, Antonio and Tenenbaum, Joshua B and Mordatch, Igor},
  booktitle={Forty-first International Conference on Machine Learning},
  year={2023}
}

@inproceedings{liang2024encouraging,
  title={Encouraging divergent thinking in large language models through multi-agent debate},
  author={Liang, Tian and He, Zhiwei and Jiao, Wenxiang and Wang, Xing and Wang, Yan and Wang, Rui and Yang, Yujiu and Shi, Shuming and Tu, Zhaopeng},
  booktitle={Proceedings of the 2024 conference on empirical methods in natural language processing},
  pages={17889--17904},
  year={2024}
}

@article{chan2023chateval,
  title={Chateval: Towards better llm-based evaluators through multi-agent debate},
  author={Chan, Chi-Min and Chen, Weize and Su, Yusheng and Yu, Jianxuan and Xue, Wei and Zhang, Shanghang and Fu, Jie and Liu, Zhiyuan},
  journal={arXiv preprint arXiv:2308.07201},
  year={2023}
}

@article{subramaniam2024debategpt,
  title={DebateGPT: Fine-tuning large language models with multi-agent debate supervision},
  author={Subramaniam, Vighnesh and Torralba, Antonio and Li, Shuang},
  year={2024}
}

@inproceedings{liu2025breaking,
  title={Breaking mental set to improve reasoning through diverse multi-agent debate},
  author={Liu, Yexiang and Cao, Jie and Li, Zekun and He, Ran and Tan, Tieniu},
  booktitle={The Thirteenth International Conference on Learning Representations},
  year={2025}
}

@article{estornell2024acc,
  title={Acc-debate: An actor-critic approach to multi-agent debate},
  author={Estornell, Andrew and Ton, Jean-Francois and Yao, Yuanshun and Liu, Yang},
  journal={arXiv e-prints},
  pages={arXiv--2411},
  year={2024}
}

@article{lewis2020retrieval,
  title={Retrieval-augmented generation for knowledge-intensive nlp tasks},
  author={Lewis, Patrick and Perez, Ethan and Piktus, Aleksandra and Petroni, Fabio and Karpukhin, Vladimir and Goyal, Naman and K{\"u}ttler, Heinrich and Lewis, Mike and Yih, Wen-tau and Rockt{\"a}schel, Tim and others},
  journal={Advances in neural information processing systems},
  volume={33},
  pages={9459--9474},
  year={2020}
}

@article{wu2024retrieval,
  title={Retrieval-augmented generation for natural language processing: A survey},
  author={Wu, Shangyu and Xiong, Ying and Cui, Yufei and Wu, Haolun and Chen, Can and Yuan, Ye and Huang, Lianming and Liu, Xue and Kuo, Tei-Wei and Guan, Nan and others},
  journal={arXiv preprint arXiv:2407.13193},
  year={2024}
}

@article{packer2023memgpt,
  title={MemGPT: Towards LLMs as Operating Systems.},
  author={Packer, Charles and Fang, Vivian and Patil, Shishir\_G and Lin, Kevin and Wooders, Sarah and Gonzalez, Joseph\_E},
  year={2023},
  journal={arXiv e-prints},
}

@inproceedings{zhong2024memorybank,
  title={Memorybank: Enhancing large language models with long-term memory},
  author={Zhong, Wanjun and Guo, Lianghong and Gao, Qiqi and Ye, He and Wang, Yanlin},
  booktitle={Proceedings of the AAAI Conference on Artificial Intelligence},
  volume={38},
  number={17},
  pages={19724--19731},
  year={2024}
}

@article{xu2025mem,
  title={A-mem: Agentic memory for llm agents},
  author={Xu, Wujiang and Mei, Kai and Gao, Hang and Tan, Juntao and Liang, Zujie and Zhang, Yongfeng},
  journal={arXiv preprint arXiv:2502.12110},
  year={2025}
}

@article{modarressi2023ret,
  title={Ret-llm: Towards a general read-write memory for large language models},
  author={Modarressi, Ali and Imani, Ayyoob and Fayyaz, Mohsen and Sch{\"u}tze, Hinrich},
  journal={arXiv preprint arXiv:2305.14322},
  year={2023}
}

@inproceedings{park2023generative,
  title={Generative agents: Interactive simulacra of human behavior},
  author={Park, Joon Sung and O'Brien, Joseph and Cai, Carrie Jun and Morris, Meredith Ringel and Liang, Percy and Bernstein, Michael S},
  booktitle={Proceedings of the 36th annual acm symposium on user interface software and technology},
  pages={1--22},
  year={2023}
}

@inproceedings{ramamonjison2023nl4opt,
  title={Nl4opt competition: Formulating optimization problems based on their natural language descriptions},
  author={Ramamonjison, Rindranirina and Yu, Timothy and Li, Raymond and Li, Haley and Carenini, Giuseppe and Ghaddar, Bissan and He, Shiqi and Mostajabdaveh, Mahdi and Banitalebi-Dehkordi, Amin and Zhou, Zirui and others},
  booktitle={NeurIPS 2022 competition track},
  pages={189--203},
  year={2023},
  organization={PMLR}
}

@article{jiang2024llmopt,
  title={LLMOPT: Learning to Define and Solve General Optimization Problems from Scratch},
  author={Jiang, Caigao and Shu, Xiang and Qian, Hong and Lu, Xingyu and Zhou, Jun and Zhou, Aimin and Yu, Yang},
  journal={arXiv preprint arXiv:2410.13213},
  year={2024}
}

@article{lu2025optmath,
  title={Optmath: A scalable bidirectional data synthesis framework for optimization modeling},
  author={Lu, Hongliang and Xie, Zhonglin and Wu, Yaoyu and Ren, Can and Chen, Yuxuan and Wen, Zaiwen},
  journal={arXiv preprint arXiv:2502.11102},
  year={2025}
}

@article{zhou2025auto,
  title={Auto-formulating dynamic programming problems with large language models},
  author={Zhou, Chenyu and Yang, Jingyuan and Xin, Linwei and Chen, Yitian and He, Ziyan and Ge, Dongdong},
  journal={arXiv preprint arXiv:2507.11737},
  year={2025}
}

@article{wu2025step,
  title={Step-Opt: Boosting Optimization Modeling in LLMs through Iterative Data Synthesis and Structured Validation},
  author={Wu, Yang and Zhang, Yifan and Wu, Yurong and Wang, Yuran and Zhang, Junkai and Cheng, Jian},
  journal={arXiv preprint arXiv:2506.17637},
  year={2025}
}

@article{chen2025solver,
  title={Solver-Informed RL: Grounding Large Language Models for Authentic Optimization Modeling},
  author={Chen, Yitian and Xia, Jingfan and Shao, Siyu and Ge, Dongdong and Ye, Yinyu},
  journal={arXiv preprint arXiv:2505.11792},
  year={2025}
}

@inproceedings{deng24cafa,
  title={CAFA: Coding as Auto-Formulation Can Boost Large Language Models in Solving Linear Programming Problem},
  author={Deng, Haoxuan and Zheng, Bohao and Jiang, Yirui and Tran, Trung Hieu},
  booktitle={The 4th Workshop on Mathematical Reasoning and AI at NeurIPS'24}
}

@article{zhou2025steporlm,
  title={Steporlm: A self-evolving framework with generative process supervision for operations research language models},
  author={Zhou, Chenyu and Xu, Tianyi and Lin, Jianghao and Ge, Dongdong},
  journal={arXiv preprint arXiv:2509.22558},
  year={2025}
}

@article{chen2025compo,
  title={ComPO: Preference alignment via comparison oracles},
  author={Chen, Peter and Chen, Xi and Yin, Wotao and Lin, Tianyi},
  journal={arXiv preprint arXiv:2505.05465},
  year={2025}
}

@article{wang2025ormind,
  title={ORMind: A Cognitive-Inspired End-to-End Reasoning Framework for Operations Research},
  author={Wang, Zhiyuan and Chen, Bokui and Huang, Yinya and Cao, Qingxing and He, Ming and Fan, Jianping and Liang, Xiaodan},
  journal={arXiv preprint arXiv:2506.01326},
  year={2025}
}

@article{kong2025alphaopt,
  title={AlphaOPT: Formulating Optimization Programs with Self-Improving LLM Experience Library},
  author={Kong, Minwei and Qu, Ao and Guo, Xiaotong and Ouyang, Wenbin and Jiang, Chonghe and Zheng, Han and Ma, Yining and Zhuang, Dingyi and Tang, Yuhan and Li, Junyi and others},
  journal={arXiv preprint arXiv:2510.18428},
  year={2025}
}

@article{wang2025mirix,
  title={Mirix: Multi-agent memory system for llm-based agents},
  author={Wang, Yu and Chen, Xi},
  journal={arXiv preprint arXiv:2507.07957},
  year={2025}
}

@article{peters2022world,
  title={UN world food programme: toward zero hunger with analytics},
  author={Peters, Koen and Silva, S{\'e}rgio and Wolter, Tim Sergio and Anjos, Luis and van Ettekoven, Nina and Combette, {\'E}ric and Melchiori, Anna and Fleuren, Hein and den Hertog, Dick and Ergun, {\"O}zlem},
  journal={Informs journal on applied analytics},
  volume={52},
  number={1},
  pages={8--26},
  year={2022},
  publisher={INFORMS}
}

@article{holland2017ups,
  title={UPS optimizes delivery routes},
  author={Holland, Chuck and Levis, Jack and Nuggehalli, Ranganath and Santilli, Bob and Winters, Jeff},
  journal={Interfaces},
  volume={47},
  number={1},
  pages={8--23},
  year={2017},
  publisher={INFORMS}
}

@article{petropoulos2024operational,
  title={Operational Research: methods and applications},
  author={Petropoulos, Fotios and Laporte, Gilbert and Aktas, Emel and Alumur, Sibel A and Archetti, Claudia and Ayhan, Hayriye and Battarra, Maria and Bennell, Julia A and Bourjolly, Jean-Marie and Boylan, John E and others},
  journal={Journal of the Operational Research Society},
  volume={75},
  number={3},
  pages={423--617},
  year={2024},
  publisher={Taylor \& Francis}
}

@article{singh2012overview,
  title={An overview of the optimization modelling applications},
  author={Singh, Ajay},
  journal={Journal of Hydrology},
  volume={466},
  pages={167--182},
  year={2012},
  publisher={Elsevier}
}

@article{wang2023voyager,
  title={Voyager: An open-ended embodied agent with large language models},
  author={Wang, Guanzhi and Xie, Yuqi and Jiang, Yunfan and Mandlekar, Ajay and Xiao, Chaowei and Zhu, Yuke and Fan, Linxi and Anandkumar, Anima},
  journal={arXiv preprint arXiv:2305.16291},
  year={2023}
}

@article{huang2024mamo,
  title={{MAMO}: A Mathematical Modeling Benchmark with Solvers},
  author={Huang, Xuhan and Shen, Qingning and Hu, Yan and Gao, Anningzhe and Wang, Benyou},
  journal={arXiv preprint},
  year={2024}
}

@inproceedings{xiao2024chainofexperts,
  author = {Ziyang Xiao and Dongxiang Zhang and Yangjun Wu and Lilin Xu and Yuan J. Wang and Xiongwei Han and Xiaojin Fu and Tao Zhong and Jia Zeng and Mingli Song and Gang Chen},
  title = {{Chain-of-Experts: When LLMs Meet Complex Operations Research Problems}},
  booktitle = {International Conference on Learning Representations (ICLR)},
  year = {2024}
}

@misc{openai2024gpt4o,
  title        = {GPT-4o System Card},
  author       = {OpenAI},
  year         = {2024},
  howpublished = {\url{https://openai.com/research/gpt-4o}},
}

@article{team2025kimi,
  title={Kimi k2: Open agentic intelligence},
  author={Team, Kimi and Bai, Yifan and Bao, Yiping and Chen, Guanduo and Chen, Jiahao and Chen, Ningxin and Chen, Ruijue and Chen, Yanru and Chen, Yuankun and Chen, Yutian and others},
  journal={arXiv preprint arXiv:2507.20534},
  year={2025}
}

@article{comanici2025gemini,
  title={Gemini 2.5: Pushing the frontier with advanced reasoning, multimodality, long context, and next generation agentic capabilities},
  author={Comanici, Gheorghe and Bieber, Eric and Schaekermann, Mike and Pasupat, Ice and Sachdeva, Noveen and Dhillon, Inderjit and Blistein, Marcel and Ram, Ori and Zhang, Dan and Rosen, Evan and others},
  journal={arXiv preprint arXiv:2507.06261},
  year={2025}
}

@article{guo2025deepseek,
  title={Deepseek-r1: Incentivizing reasoning capability in llms via reinforcement learning},
  author={Guo, Daya and Yang, Dejian and Zhang, Haowei and Song, Junxiao and Zhang, Ruoyu and Xu, Runxin and Zhu, Qihao and Ma, Shirong and Wang, Peiyi and Bi, Xiao and others},
  journal={arXiv preprint arXiv:2501.12948},
  year={2025}
}

@misc{deepseekai2025deepseekv3technicalreport,
      title={DeepSeek-V3 Technical Report}, 
      author={DeepSeek-AI and Aixin Liu and Bei Feng and Bing Xue and Bingxuan Wang and Bochao Wu and Chengda Lu and Chenggang Zhao and Chengqi Deng and Chenyu Zhang and Chong Ruan and Damai Dai and Daya Guo and Dejian Yang and Deli Chen and Dongjie Ji and Erhang Li and Fangyun Lin and Fucong Dai and Fuli Luo and Guangbo Hao and Guanting Chen and Guowei Li and H. Zhang and Han Bao and Hanwei Xu and Haocheng Wang and Haowei Zhang and Honghui Ding and Huajian Xin and Huazuo Gao and Hui Li and Hui Qu and J. L. Cai and Jian Liang and Jianzhong Guo and Jiaqi Ni and Jiashi Li and Jiawei Wang and Jin Chen and Jingchang Chen and Jingyang Yuan and Junjie Qiu and Junlong Li and Junxiao Song and Kai Dong and Kai Hu and Kaige Gao and Kang Guan and Kexin Huang and Kuai Yu and Lean Wang and Lecong Zhang and Lei Xu and Leyi Xia and Liang Zhao and Litong Wang and Liyue Zhang and Meng Li and Miaojun Wang and Mingchuan Zhang and Minghua Zhang and Minghui Tang and Mingming Li and Ning Tian and Panpan Huang and Peiyi Wang and Peng Zhang and Qiancheng Wang and Qihao Zhu and Qinyu Chen and Qiushi Du and R. J. Chen and R. L. Jin and Ruiqi Ge and Ruisong Zhang and Ruizhe Pan and Runji Wang and Runxin Xu and Ruoyu Zhang and Ruyi Chen and S. S. Li and Shanghao Lu and Shangyan Zhou and Shanhuang Chen and Shaoqing Wu and Shengfeng Ye and Shengfeng Ye and Shirong Ma and Shiyu Wang and Shuang Zhou and Shuiping Yu and Shunfeng Zhou and Shuting Pan and T. Wang and Tao Yun and Tian Pei and Tianyu Sun and W. L. Xiao and Wangding Zeng and Wanjia Zhao and Wei An and Wen Liu and Wenfeng Liang and Wenjun Gao and Wenqin Yu and Wentao Zhang and X. Q. Li and Xiangyue Jin and Xianzu Wang and Xiao Bi and Xiaodong Liu and Xiaohan Wang and Xiaojin Shen and Xiaokang Chen and Xiaokang Zhang and Xiaosha Chen and Xiaotao Nie and Xiaowen Sun and Xiaoxiang Wang and Xin Cheng and Xin Liu and Xin Xie and Xingchao Liu and Xingkai Yu and Xinnan Song and Xinxia Shan and Xinyi Zhou and Xinyu Yang and Xinyuan Li and Xuecheng Su and Xuheng Lin and Y. K. Li and Y. Q. Wang and Y. X. Wei and Y. X. Zhu and Yang Zhang and Yanhong Xu and Yanhong Xu and Yanping Huang and Yao Li and Yao Zhao and Yaofeng Sun and Yaohui Li and Yaohui Wang and Yi Yu and Yi Zheng and Yichao Zhang and Yifan Shi and Yiliang Xiong and Ying He and Ying Tang and Yishi Piao and Yisong Wang and Yixuan Tan and Yiyang Ma and Yiyuan Liu and Yongqiang Guo and Yu Wu and Yuan Ou and Yuchen Zhu and Yuduan Wang and Yue Gong and Yuheng Zou and Yujia He and Yukun Zha and Yunfan Xiong and Yunxian Ma and Yuting Yan and Yuxiang Luo and Yuxiang You and Yuxuan Liu and Yuyang Zhou and Z. F. Wu and Z. Z. Ren and Zehui Ren and Zhangli Sha and Zhe Fu and Zhean Xu and Zhen Huang and Zhen Zhang and Zhenda Xie and Zhengyan Zhang and Zhewen Hao and Zhibin Gou and Zhicheng Ma and Zhigang Yan and Zhihong Shao and Zhipeng Xu and Zhiyu Wu and Zhongyu Zhang and Zhuoshu Li and Zihui Gu and Zijia Zhu and Zijun Liu and Zilin Li and Ziwei Xie and Ziyang Song and Ziyi Gao and Zizheng Pan},
      year={2025},
      eprint={2412.19437},
      archivePrefix={arXiv},
      primaryClass={cs.CL},
      url={https://arxiv.org/abs/2412.19437}, 
}

@misc{chen2026optenginebenchmarkinglimitsllms,
      title={OPT-Engine: Benchmarking the Limits of LLMs in Optimization Modeling via Complexity Scaling}, 
      author={Yitian Chen and Cheng Cheng and Yinan Sun and Zi Ling and Dongdong Ge},
      year={2026},
      eprint={2601.19924},
      archivePrefix={arXiv},
      primaryClass={cs.CL},
      url={https://arxiv.org/abs/2601.19924}, 
}
